\documentclass{amsart}
\usepackage{amssymb,amsfonts}
\usepackage{enumerate}
\usepackage{enumitem}
\usepackage{mathrsfs}
\usepackage{graphicx}
\usepackage{mathtools}
\usepackage{leftidx}
\usepackage{tikz-cd}
\usepackage{amsmath,calligra,mathrsfs}
\usepackage{bm}
\usepackage{url}
\newtheorem{thm}{Theorem}[section]
\newtheorem{cor}[thm]{Corollary}
\newtheorem{prop}[thm]{Proposition}
\newtheorem{lem}[thm]{Lemma}

\theoremstyle{definition}
\newtheorem{defn}[thm]{Definition}

\newtheorem{exmp}[thm]{Example}

\theoremstyle{remark}
\newtheorem{rem}[thm]{Remark}

\makeatletter
\let\c@equation\c@thm
\makeatother
\numberwithin{equation}{section}

\def\bthm{\begin{thm}}
\def\ethm{\end{thm}}
\def\blm{\begin{lem}}
\def\elm{\end{lem}}
\def\bdf{\begin{defn}}
\def\edf{\end{defn}}
\def\bpf{\begin{proof}}
\def\epf{\end{proof}}
\def\bpp{\begin{prop}}
\def\epp{\end{prop}}
\def\bcor{\begin{cor}}
\def\ecor{\end{cor}}
\def\brm{\begin{rem}}
\def\erm{\end{rem}}
\def\beg{\begin{exmp}}
\def\eeg{\end{exmp}}
\def\bD{\mathbb{D}}

\def\bG{\mathbb{G}}

\def\bN{\mathbb{N}}

\def\bQ{\mathbb{Q}}

\def\bZ{\mathbb{Z}}
\def\cA{\mathcal{A}}
\def\cB{\mathcal{B}}
\def\cC{\mathcal{C}}
\def\cD{\mathcal{D}}
\def\cE{\mathcal{E}}
\def\cF{\mathcal{F}}

\def\cL{\mathcal{L}}
\def\cM{\mathcal{M}}
\def\cN{\mathcal{N}}
\def\cO{\mathcal{O}}
\def\cP{\mathcal{P}}

\def\cR{\mathcal{R}}

\def\cT{\mathcal{T}}

\def\scI{\mathscr{I}}

\def\scL{\mathscr{L}}

\newcommand{\raq}{\,\rightarrow \,}

\newcommand{\xraq}[2][]{\, \xrightarrow[#1]{#2} \,}

\newcommand{\ra}{\rightarrow}

\newcommand{\rinto}{\hookrightarrow}

\newcommand{\ronto}{\twoheadrightarrow}

\newcommand{\rsa}{\stackrel{\sim}{\rightarrow}}

\newcommand{\xra}[2][]{\xrightarrow[#1]{#2}}

\newcommand{\ie}{{\it i.e.}}
\newcommand{\eg}{{\it e.g.}}

\newcommand{\Mod}{{\rm Mod}}

\newcommand{\Ab}{{\rm Ab}}

\newcommand{\Ch}{{\rm Ch}}
\newcommand{\cone}{{\rm cone}}

\newcommand{\Tor}{{\rm Tor}}

\newcommand{\Extcom}{\underline{{\rm Ext}}}

\newcommand{\Ob}{{\rm Ob}}
\newcommand{\op}{{\rm op}}

\newcommand{\id}{{\rm id}}

\newcommand{\Hom}{{\rm Hom}}
\newcommand{\Homcom}{\underline{{\rm Hom}}}

\newcommand{\cHom}{\mathscr{H}\text{\kern -3pt {\calligra\large om}}\,}

\newcommand{\RHom}{{\bm R}{\rm Hom}}
\newcommand{\RHomcom}{{\bm R} \underline{{\rm Hom}}}

\newcommand{\Spec}{{\rm Spec}}

\newcommand{\Proj}{{\rm Proj}}
\newcommand{\PProj}{\mathpzc{Proj}}

\newcommand{\QCoh}{{\rm QCoh}}

\newcommand{\Dsuit}{\cD^{\tiny \mbox{$\spadesuit $}}}

\newcommand{\Dbcoh}{\cD^b_{{\rm coh}}}

\newcommand{\Dmcoh}{\cD^{-}_{{\rm coh}}}

\newcommand{\Dpc}{\cD_{{\rm pc}}}

\newcommand{\Dperf}{\cD_{{\rm perf}}}

\newcommand{\perf}{{\rm perf}}
\newcommand{\coh}{{\rm coh}}

\newcommand{\NS}{{\rm NS}}

\newcommand{\Gr}{{\rm Gr}}
\newcommand{\GrA}{{\rm Gr}(A)}
\newcommand{\gr}{{\rm gr}}

\newcommand{\QpGr}{{\rm Q}^+{\rm Gr}}

\newcommand{\qpgr}{{\rm q}^+{\rm gr}}

\newcommand{\Torp}{{\rm Tor}^+}

\newcommand{\IIT}{I^{\infty}\text{-}{\rm Tor}}

\newcommand{\cIIT}{\scI^{\infty}\text{-}{\rm Tor}}

\newcommand{\IpIT}{(I^+)^{\infty}\text{-}{\rm Tor}}

\newcommand{\Imcomp}{I^-\text{-}{\rm comp}}

\newcommand{\RGam}{{\bm R}\Gamma}

\newcommand{\RGI}{{\bm R}\Gamma_{I}}

\newcommand{\Ce}{\check{\cC}}
\newcommand{\CI}{\check{\cC}_I}

\newcommand{\LLam}{{\bm L}\Lambda}

\newcommand{\ITR}{I\text{-}{\rm triv}}
\newcommand{\IpTR}{I^+\text{-}{\rm triv}}
\newcommand{\ImTR}{I^-\text{-}{\rm triv}}

\newcommand{\cITR}{\scI\text{-}{\rm triv}}
\newcommand{\cIpTR}{\scI^+\text{-}{\rm triv}}

\newcommand{\Fga}{\cF_{\geq a}}
\newcommand{\Fla}{\cF_{\leq a}}
\newcommand{\Flma}{\cF_{\leq -a}}

\newcommand{\Fgw}{\cF_{\geq w}}

\newcommand{\Flmw}{\cF_{\leq -w}}

\newcommand{\DAO}{\bD_A}

\DeclareMathAlphabet{\mathpzc}{OT1}{pzc}{m}{it}

\title{Weight truncation for wall-crossings in birational cobordisms}
\author{Wai-Kit Yeung}
\address{Department of Mathematics,
	Indiana University, Bloomington, IN 47405, USA}
\email{yeungw@iu.edu}
\begin{document}

\begin{abstract} 
We develop the technique of weight truncation in the context of wall-crossings in birational cobordisms, parallel to that in \cite{HL15, BFK19}. More precisely, for each such wall-crossing, we embed the bounded above derived category of coherent sheaves of the semistable part as a semi-orthogonal summand of that of the stack in question. Our construction does not require any smoothness assumptions, and exhibits a strong symmetry across the two sides of the wall-crossing. As an application, we show that for wall-crossings satisfying suitable regularity conditions, a certain duality of local cohomology complexes implies the existence of a fully faithful functor/equivalence between the derived categories under wall-crossings.
\end{abstract}

\maketitle


\tableofcontents

\section{Introduction}

A central problem in the geometry of derived categories is to study the change in the derived categories under birational transforms. In recent years, significant progress on this problem have been obtained in the setting of variations of GIT quotients, or more broadly wall-crossings in moduli problems. 
This is a convenient setting because, on the one hand, the spaces of semistable objects are often birational to each other under wall-crossings; on the other hand, the fact that they arise as different open subspaces of the same space $\mathfrak{X}$ allows one to relate their derived categories via their mutual relations to the derived category of the space $\mathfrak{X}$.

For example, if $G$ is a reductive group acting on a quasi-projective variety $X$, then each choice of $L \in \NS^G(X)_{\bQ}$ determines an open substack $\mathfrak{X}^{ss} = \mathfrak{X}^{ss}(L)$ of the stack $\mathfrak{X} := [X / G]$.
By a standard abuse of language, we will refer to $\mathfrak{X}^{ss}$ as the GIT quotient, the usual GIT quotient being its scheme-theoretic categorical quotient.
As an open substack, the derived category of $\mathfrak{X}^{ss}$ is a localization of that of $\mathfrak{X}$. 
As $L$ changes, say from $L^-$ to $L^+$, then the derived categories of $\mathfrak{X}^{\pm} := \mathfrak{X}^{ss}(L^{\pm})$ are thus two different localizations of the same derived category. This setting is a baseline for the comparison of the derived categories of $\mathfrak{X}^-$ and $\mathfrak{X}^+$.

Specifically, one may compare the derived categories of $\mathfrak{X}^-$ and $\mathfrak{X}^+$ by finding suitable full triangulated subcategories $\cE^{\pm} \subset \Dbcoh(\mathfrak{X})$ such that the restriction functors $(j^{\pm})^* : \cE^{\pm} \ra \Dbcoh(\mathfrak{X}^{\pm})$ are  equivalences. There is often a direct comparison between $\cE^-$ and $\cE^+$, which in turn allows one to compare $\Dbcoh(\mathfrak{X}^-)$ and $\Dbcoh(\mathfrak{X}^+)$.

Weight truncation, also known as grade restriction rule or GIT window (see, \eg, \cite{HL15, BFK19}), is a particular way to realize this idea. In this paper, we develop the technique of weight truncation, focusing on the abelian case $G = \bG_m$, or more precisely the case of birational cobordism, where the element $L \in \NS^{\bG_m}(X)_{\bQ}$ changes via twisting by a (fractional) character of $\bG_m$.
Indeed, there are several techniques to reduce the more general case to this abelian case. For example, the ``master space construction'' of Thaddeus \cite{Tha96} realizes every variation of GIT quotients as a birational cobordism. On the other hand, the structure of HKKN stratification essentially allows for a strata-by-strata reduction to the abelian case. This last picture also extends beyond GIT quotients (see, \eg, \cite{HL}). As such, the abelian case has a special significance.

Compared with \cite{HL15, BFK19}, our approach has several advantages. Firstly, our construction of weight truncation does not require any smoothness assumption, although some regularity condition is still required in order to guarantee certain desired properties. Secondly, our construction exhibits a manifest symmetry across the wall. Thirdly, our construction carries over to the noncommutative case, in the spririt of noncommutative projective geometry (see, \eg, \cite{AZ94, Orl09}), although this will not be made explicit in this paper.

Let us now explain our construction of weight truncation. 
Given a $\bG_m$-space $X$, consider a family $L(t) \in \NS^{\bG_m}(X)_{\bQ}$ obtained by twisted a given $L$ by fractional characters $t \in \bQ$. Under a wall-crossing from $t^-$ to $t^+$, passing through the wall $t_0$, the semistable loci satisfy $X^{ss}(L(t^-)) \subset X^{ss}(L(t_0)) \supset X^{ss}(L(t^+))$. Denote by $W := X^{ss}(L(t_0))$, and $Y := W // \bG_m$ the scheme-theoretic categorical quotient. Then we have $W = \Spec_Y(\cA)$ for a sheaf of $\bZ$-graded rings on $Y$, and the (stacky) GIT quotients across the wall are
\begin{equation}  \label{birat_cobord_log_flip_diag_intro}
\left( 
\begin{tikzcd} [row sep = 0, column sep = 0]
{\rm [} X^{ss}(L(t^-)) \, / \, \bG_m {\rm ]} \ar[rd]
& &  {\rm [} X^{ss}(L(t^+)) \, / \, \bG_m {\rm ]} \ar[ld] \\
& Y &
\end{tikzcd}
\right)
\, \cong \, 
\left( 
\begin{tikzcd} [row sep = 0, column sep = 5]
\PProj^-_Y(\cA) \ar[rd]
& &  \PProj^+_Y(\cA) \ar[ld] \\
& Y &
\end{tikzcd}
\right)
\end{equation}
where $\PProj^{+}_Y(\cA)$ is a stacky version of the scheme-theoretic GIT quotient $\Proj^+_Y(\cA) := \Proj_Y(\cA_{\geq 0})$, and similarly for  $\PProj^{-}_Y(\cA)$. We will write $\mathfrak{X}^{\pm} := \PProj^{\pm}_Y(\cA)$.

Thus, the study of wall-crossing in this abelian case is reduced to the study of a sheaf of $\bZ$-graded rings. Since all our constructions are local on $Y$, we may assume that $Y$ is affine, so that we may focus on a Noetherian $\bZ$-graded ring $A$ ({\it cf.} Proposition \ref{Noeth_gr_ring} below for the Noetherian property).
The category of quasi-coherent sheaves on $\mathfrak{X} := [\Spec \, A / \bG_m]$ in this case is simply the category $\GrA$ of graded $A$-modules. 
Thus, we are reduced to studying the derived category $\cD(\GrA)$ of graded $A$-modules.


Let $I^+ := A_{>0} \cdot A$ and $I^- := A_{<0} \cdot A$. Then the semistable loci $X^{ss}(L^{\pm})$ are the complement of the $\bG_m$-invariant closed subsets defined by the graded ideals $I^{\pm}$. Thus, if we denote by $j^{\pm} : \mathfrak{X}^{\pm} \subset \mathfrak{X}$ the open inclusions, then for any $M \in \cD(\GrA) \simeq \cD(\QCoh(\mathfrak{X}))$, the derived pushforward ${\bm R}j^{+}_*(j^{+})^*(M)$ is computed by a certain \v{C}ech complex $\Ce_{I^+}(M)$ (see Definition \ref{Cech_cplx_def} below), which sits in an exact triangle
\begin{equation} \label{RGI_CI_seq_intro}
\ldots \raq \RGam_{I^+}(M) \xraq{\epsilon} M \xraq{\eta} \Ce_{I^+}(M) \xraq{\delta} \RGam_{I^+}(M)[1] \raq \ldots 
\end{equation}

In fact, this triangle is the decomposition triangle for the semi-orthogonal decomposition
\begin{equation}  \label{torsion_SOD_intro}
\cD(\GrA) \, = \, \langle \, \cD_{\IpTR}(\GrA) \, , \,  \cD_{\Torp}(\GrA) \, \rangle
\end{equation}
where $\cD_{\Torp}(\GrA)$ are objects of $\cD(\QCoh(\mathfrak{X}))$ supported on the unstable locus, and 
$\cD_{\IpTR}(\GrA)$ is equivalent to  $\cD(\QCoh(\mathfrak{X}^+))$ via the functor ${\bm R}j^{\pm}_*$.

While this gives an embedding of $\cD(\QCoh(\mathfrak{X}^+))$ as a semi-orthogonal summand of $\cD(\GrA)$, this embedding is usually very far%
\footnote{By the exact triangle \eqref{RGI_CI_seq_intro}, this is the same as asking whether the local cohomology modules $H^j(\RGam_{I^+}(M))$ are finitely generated. For a result along this line, see, \eg, \cite[Tag 0BJV]{Sta}.} from preserving $\Dbcoh$:
\begin{equation}  \label{Cech_coherent_intro_1}
\parbox{40em}{For $M \in \Dbcoh(\GrA)$, the \v{C}ech complex $\Ce_{I^+}(M)$ is almost never in $\Dbcoh(\GrA)$.}
\end{equation}

The \v{C}ech complex $\Ce_{I^+}(M)$ also plays a role for the classical GIT quotient. Indeed, let $X^+ := \Proj(A_{\geq 0})$, and $\pi^+ : X^+ \ra Y = \Spec \, A_0$ the projection, then since the cohomology of quasi-coherent sheaves can be computed by the \v{C}ech complex, we have
\begin{equation}  \label{CI_Rpi_intro}
{\bm R}\pi^+_* (\widetilde{M(i)}) \, \cong \, \Ce_{I^+}(M) _i
\end{equation}

It is a classical fact that the derived pushforward of the projective morphism $\pi^+ : X^+ \ra Y$ sends $\Dbcoh(X^+)$ to $\Dbcoh(Y)$. Thus, by \eqref{CI_Rpi_intro}, each weight component $\Ce_{I^+}(M)_i$ of the \v{C}ech complex is in $\Dbcoh(A_0)$. In fact, more is true (see the proof of Proposition \ref{weight_Cech_Dbcoh} below):
\begin{equation}  \label{Cech_coherent_intro_2}
\parbox{40em}{If $M \in \Dbcoh(\GrA)$, then for each $w \in \bZ$, take the restriction $\Ce_{I^+}(M)_{\geq w}$ of the \v{C}ech complex to weights $i \geq w$, and consider it as an object in $\cD(\Gr(A_{\geq 0}))$. Then we have $\Ce_{I^+}(M)_{\geq w} \in \Dbcoh(\Gr(A_{\geq 0}))$.}
\end{equation}

A comparison of \eqref{Cech_coherent_intro_1} and \eqref{Cech_coherent_intro_2} suggests that, by considering the restriction $\Ce_{I^+}(M)_{\geq w}$ of the \v{C}ech complex, one might be able to obtain an embedding of $\cD(\QCoh(\mathfrak{X}^+))$ into $\cD(\GrA)$ that has more chance of preserving $\Dbcoh(-)$. This is achieved by the technique of weight truncation. For this, we will reformulate the observation \eqref{Cech_coherent_intro_2}.

In \eqref{Cech_coherent_intro_2}, we considered the restriction $\Ce_{I^+}(M)_{\geq w}$ as an object in $\cD(\Gr(A_{\geq 0}))$. However, the passage from $A$ to $A_{\geq 0}$ is somewhat unnatural, and it is difficult to compare their derived categories. In fact, there is a better way to organize the data of the complexes $\Ce_{I^+}(M)_{i}$ for $i \geq w$. For this, we use the language of a small pre-additive category. 

By definition, a pre-additive category is a category whose Hom sets have structures of abelian groups, and whose composition maps are bilinear. As in \cite{Mit72}, one may view a small pre-additive category as an associative algebra with many objects, and as such, the usual notions about associative algebras, such as left/right modules, left/right Noetherian properites, tensor product of modules, derived categories of modules, etc, carries over to the case of small pre-additive categories. See Appendix \ref{app_mod_preadd} for a summary of all the results needed for our discussion.

Let $\cF$ be the pre-additive category with objects set $\Ob(\cF) = \bZ$, and with Hom sets $\cF(i,j) := A_{i-j}$. Compositions are defined by the multiplication in $A$. It is clear that a (right) module of $\cF$ is the same as a graded module over $A$. Let $\Fgw \subset \cF$ be the full subcategory with object set $\bZ_{\geq w}$. Then for any graded module $M$, the data $M_{\geq w}$ is naturally organized into a right module over $\Fgw$. The observation \eqref{Cech_coherent_intro_2} can then be rewritten in the following form:
\begin{equation}  \label{Cech_coherent_intro_3}
\parbox{40em}{If $M \in \Dbcoh(\GrA)$, then for each $w \in \bZ$, the restriction $\Ce_{I^+}(M)_{\geq w}$ is in $\Dbcoh(\Fgw)$.}
\end{equation}

This observation is the crucial entry into the proof of the following result, which asserts that the ``weight truncated'' analogue of \eqref{torsion_SOD_intro} always preserve $\Dbcoh$.
\bthm [= Corollary \ref{weight_Cech_Dsuitcoh_SOD}, Theorem \ref{I_triv_weight_equiv_Dsuitcoh}] \label{Dbcoh_Fgw_SOD_intro_thm}
There is a semi-orthogonal decomposition
\begin{equation}  \label{Dbcoh_Fgw_SOD_intro}
\Dbcoh(\Fgw) = \langle \, \cD^b_{{\rm coh}, \, \IpTR}(\Fgw) \, , \, \cD^b_{{\rm coh}, \, \Torp}(\Fgw) \, \rangle
\end{equation}
Moreover, there is an exact equivalence $\cD^b_{{\rm coh}, \, \IpTR}(\Fgw) \simeq \Dbcoh(\mathfrak{X}^+)$. 
\ethm


One can then use Theorem \ref{Dbcoh_Fgw_SOD_intro_thm} to obtain results on $\cD(\GrA)$. Indeed, since $\Fgw$ is a full subcategory of $\cF$, one may use (left) Kan extensions to obtain a fully faithful functor $\scL_{[\geq w]} : \cD(\Fgw) \rinto \cD(\cF) \simeq \cD(\GrA)$. By proving a certain Noetherian property of $\Fgw$ (see Corollary \ref{Fgw_Noeth} below), one can show that this embedding always preserves $\Dmcoh$. Hence, a version of Theorem \ref{Dbcoh_Fgw_SOD_intro_thm} with $\Dbcoh$ replaced by $\Dmcoh$ then implies that there is a three-term semi-orthogonal decomposition (see Theorem \ref{three_term_SOD_Dmcoh})
\begin{equation*}
\Dmcoh(\GrA) \, = \, \langle \, \cD^-_{\coh,<w}(\GrA) \, , \, \mathscr{L}_{[\geq w]} ( \cD^-_{\coh, \, \IpTR}(\Fgw) )  \, , \,  \mathscr{L}_{[\geq w]} ( \cD^-_{\coh, \, \Torp}(\Fgw) ) \, \rangle
\end{equation*}
with the middle component equivalent to $\Dmcoh(\mathfrak{X}^+)$.
Similar result was obtained in \cite{HL15b} using different methods. We believe that these give rise to the same semi-orthogonal decomposition. We also give a sufficient condition for this semi-orthogonal decomposition to restrict to $\Dbcoh(\GrA)$ (see Theorem \ref{cond_then_reg}). In particular, this weakens the ``Assumption (A)'' in \cite{HL15}. When these conditions are satisfied, we also show that the semi-orthogonal decomposition thus obtained coincide with the one in \cite{HL15}.

As we mentioned above, the technique of weight truncation is designed to faciliate the comparison of derived categories under wall-crossings. Theorem \ref{Dbcoh_Fgw_SOD_intro_thm} is particularly convenient for that purpose, as it exhibits a strong symmetry between the two sides of the wall-crossing. Indeed, the analogue of Theorem \ref{Dbcoh_Fgw_SOD_intro_thm} for the negative direction concerns the category $\Flmw$, which is in fact isomorphic to the opposite category of $\Fgw$. 
In other words, $\Dbcoh(\mathfrak{X}^+)$ is a semi-orthogonal summand of the derived category $\Dbcoh(\Fgw)$ of \emph{right} $\Fgw$-modules; while $\Dbcoh(\mathfrak{X}^-)$ is a semi-orthogonal summand of the derived category $\Dbcoh((\Fgw)^{\op})$ of \emph{left} $\Fgw$-modules.
This suggests one to relate the two via a duality functor
\begin{equation}  \label{bD_Fgw_intro}
\bD_{\Fgw}  :  \cD(\Fgw)^{\op} \ra \cD((\Fgw)^{\op}) \simeq \cD(\Flmw) \, , \qquad \quad 
\cM \, \mapsto \, \RHom_{\Fgw}(\cM, \Fgw)
\end{equation}

We will show below that this works well when there is a certain duality between the local cohomology complexes $\RGam_{I^+}(A)$ and $\RGam_{I^-}(A)$. Let $\omega_Y^{\bullet} \in \Dbcoh(A_0)$ be a dualizing complex, and take the weight degreewise dualizing functor
\begin{equation*}
\bD_Y : \cD(\GrA)^{\op} \ra \cD(\GrA) \, , \qquad \quad \bD_Y(M)_i \simeq \RHom_{A_0}(M_{-i}, \omega_Y^{\bullet})
\end{equation*} 
Then we assume that
\begin{equation}  \label{local_cohom_duality_intro}
\parbox{40em}{There is an isomorphism $\RGam_{I^+}(A)(a)[1] \xra{\cong }\bD_Y(\RGam_{I^-}(A))$ in $\cD(\GrA)$, where $a \geq 0$.}
\end{equation}

Indeed, one of the main results of \cite{Yeu20a} is that, for a large class of flips and flops, the assumption \eqref{local_cohom_duality_intro} is satisfied for the sheaf of graded rings that controls the flip/flop, where $a=0$ for flops and $a=1$ for flips.
We have the following (see Proposition \ref{D_Fgw_ITR} below)

\bpp
If \eqref{local_cohom_duality_intro} holds for $a \geq 0$, then $\bD_{\Fgw}$ sends $\cD^-_{{\rm coh}, \, \IpTR}(\Fgw)$ to $\cD^+_{{\rm coh}, \, \ImTR}(\Flmw)$. 

If \eqref{local_cohom_duality_intro} holds for $a = 0$, then $\bD_{\Flmw}$ also sends $\cD^-_{{\rm coh}, \, \ImTR}(\Flmw)$ to $\cD^+_{{\rm coh}, \, \IpTR}(\Fgw)$. 
\epp

Thus, if the duality functor \eqref{bD_Fgw_intro} is suitably involutive, say if $\Fgw$ is ``Gorenstein'' in a suitable sense, then this result says that the derived categories behave in the expected way: it shrinks under flips and remain equivalent under flops. 
In any case, the duality functor is always involutive on $\Dperf(\Fgw)$, and we have the following (see Corollary \ref{flip_flop_Dperf})

\bcor  \label{Dperf_X_pm_intro}
Assume that the semi-orthogonal decomopsition \eqref{Dbcoh_Fgw_SOD_intro} and its negative version restrict to $\Dperf(\Fgw)$ and $\Dperf(\Flmw)$ respectively. If \eqref{local_cohom_duality_intro} holds for $a \geq 0$, then there is a fully faithful exact functor $\Dperf(\mathfrak{X}^+) \rinto \Dperf(\mathfrak{X}^-)$; if \eqref{local_cohom_duality_intro} holds for $a = 0$, then there is an exact equivalence $\Dperf(\mathfrak{X}^+) \simeq \Dperf(\mathfrak{X}^-)$.
\ecor

In fact, if $A$ is smooth, then we have $\Dbcoh(\Fgw) = \Dperf(\Fgw)$, so that the assumption of this Corollary is indeed satisfied in good cases. In general, a further sharpening of Corollary \ref{Dperf_X_pm_intro} shows that the analogous conclusions at the level of Ind-completions hold provided that certain spectral sequences converge (see Remark \ref{pairing_remark}). Such a convergence issue seems to be similar to those arising in Koszul duality. It seems possible that a formal modification of our arguments might lead to a more satisfactory statement. This might open up a way to tackle a conjecture of Bondal and Orlov (see also \cite{Yeu20a}).


\section{Derived categories of graded modules}  \label{DGrA_sec}

In this section, we recall some basic results about graded rings, in order to establish notations and conventions. A more detailed discussion may be found in \cite{Yeu20c}.

\bdf
A \emph{$\bZ$-graded ring} is a commutative ring $A$ with a $\bZ$-grading $A = \bigoplus_{n \in \bZ} A_n$. 
Here, by commutative we mean $xy = yx$, not $xy = (-1)^{|x||y|} yx$. 

A \emph{graded module} over $A$ will always mean a $\bZ$-graded module $M = \bigoplus_{n \in \bZ} M_n$.
\edf

We first recall the following result (see, e.g., \cite[Theorem 1.5.5]{BH93}):
\bpp  \label{Noeth_gr_ring}
Let $A$ be a $\bZ$-graded ring. Then the followings are equivalent:
\begin{enumerate}
	\item $A$ is a Noetherian ring; 
	\item every graded ideal of $A$ is finitely generated;
	\item $A_0$ is Noetherian, and both $A_{\geq 0}$ and $A_{\leq 0}$ are finitely generated over $A_0$;
	\item $A_0$ is Noetherian, and $A$ is finitely generated over $A_0$.
\end{enumerate}
\epp

Denote by $\Gr(A)$ the category of graded modules over $A$, whose morphisms are maps of graded modules of degree $0$.
Given two graded modules $M,N \in \Gr(A)$, then the $A$-module $M \otimes_A N$ has a natural grading where $\deg(x\otimes y) = \deg(x) + \deg(y)$ for homogeneous $x,y \in A$. Moreover, one can define a graded $A$-module $\Homcom_A(M,N)$ whose degree $i$ part is the set of $A$-linear homomorphism from $M$ to $N$ of homogeneous degree $i$. 
Thus, in particular, we have $\Hom_A(M,N) := \Hom_{\Gr(A)}(M,N) = \Homcom_A(M,N)_0$. These form the internal Hom objects with respect to the graded tensor product. 

The abelian category $\GrA$ is a Grothendieck category, with a set $\{A(-i)\}_{i \in \bZ}$ of generators. The same set is also a set of compact generators in the derived category $\cD(\GrA)$.
Since $\GrA$ is a Grothendieck category, the category of complexes has enough K-injectives (see, e.g., \cite[Tag 079P]{Sta}). Moreover, as in the ungraded case, it also has enough K-projectives (see, e.g., \cite[Tag 06XX]{Sta}).
As a result, the bifunctors $-\otimes_A - $ and $\Homcom_A(-,-)$ admit derived functors
\begin{equation*}
\begin{split}
- \otimes_{A}^{{\bm L}} - \, &: \, \cD(\GrA) \, \times \, \cD( \GrA) \raq \cD(\GrA) \\
\RHomcom_{A}(-,-) \, &: \, \cD(\GrA)^{\op} \, \times \, \cD( \GrA ) \raq \cD(\GrA)
\end{split}
\end{equation*}
which can in turn be used to define $\Extcom^{\bullet}_A(M,N)$ and $\Tor^A_{\bullet}(M,N)$.

%

We now extend some standard results on the derived categories of modules to the graded case. We start with the following

\bdf
An object $M \in \cD(\GrA)$ is said to be \emph{pseudo-coherent} if it can be represented by a bounded above complex of finitely generated projective graded $A$-modules. 
Denote by $\Dpc(\GrA) \subset \cD(\GrA)$ the full subcategory consisting of pseudo-coherent objects.
\edf

\bdf
Given a Noetherian $\bZ$-graded ring $A$, then for $\spadesuit \in \{\, \, \,   ,+,-,b\}$, define $\Dsuit_{\coh}(\GrA)$ the full subcategory of $\Dsuit(\GrA)$ consisting of complexes $M \in \Dsuit(\GrA)$ such that $H^p(M)$ is finitely generated for all $p \in \bZ$.
\edf

Then we have the following standard result (see Lemma \ref{Dbcoh_pc_lem}):

\bpp  
Given a Noetherian $\bZ$-graded ring $A$, then for any $M \in \cD(\cA)$, the followings are equivalent:
\begin{enumerate}
	\item $M \in \Dpc(\GrA)$;
	\item $M \in \Dmcoh(\GrA)$;
	\item $M$ is quasi-isomorphic to a bounded above complex of free graded modules of finite rank.
\end{enumerate}
\epp

Let $\cD_{\perf}(\GrA)$ be the smallest split-closed triangulated subcategories containing the set $\{A(-i)\}_{i \in \bZ}$ of objects. 
By \cite[Theorem 4.22]{Rou08} (see also \cite[Lemma 2.2]{Nee92}),  $\cD_{\perf}(\GrA)$ is precisely the full subcategory of compact objects in $\cD(\GrA)$. 
%

We mention the following graded analogue of 
\cite[Tag 0ATK]{Sta}, whose proof is completely parallel to the ungraded case:



\bpp  \label{tensor_in_Hom_target}
For any $N \in \Dpc(\GrA)$, $L \in \cD^+(\GrA)$, and $M \in \cD(\GrA)$ of finite Tor dimension, the canonical map
\begin{equation*}
M \otimes_A^{{\bm L}} \RHomcom_A(N,L) \raq \RHomcom_A(N, M\otimes_A^{{\bm L}} L)
\end{equation*}
is an isomorphism in $\cD(\GrA)$.
\epp

Now we briefly discuss local cohomology on a graded ring.

\bdf  \label{I_infty_torsion_def}
Let $I$ be a graded ideal in a $\bZ$-graded ring $A$. Given any graded module $M$ over $A$, an element $x\in M$ is said to be \emph{$I^{\infty}$-torsion} 
if for every $f \in I$ there exists some $n > 0$ such that $f^n x = 0$. If $I$ is finitely generated, this is equivalent to $I^n x = 0$ for some $n > 0$. The graded module $M$ is said to be \emph{$I^{\infty}$-torsion} if every element in it is $I^{\infty}$-torsion.
Denote by $\IIT \subset \Gr(A)$ the full subcategory consisting of $I^{\infty}$-torsion modules.
\edf

It is clear that $\IIT \subset \GrA$ is a Serre subcategory. Thus the full subcategory $\cD_{\IIT}(\GrA) \subset \cD(\GrA)$ is a triangulated subcategory.
If $I$ is finitely generated, then this inclusion has a right adjoint, which has a simple and useful description.
To this end, we recall the following

\bdf  \label{RGam_f_def}
Let $f_1,\ldots,f_r$ be homogeneous elements in $A$, of degrees $d_1,\ldots,d_r$ respectively. 
For any graded module $M \in \Gr(A)$, 
we define the \emph{local cohomology complex} (or \emph{extended \v{C}ech complex}) of $M$ with respect to the tuple $(f_1,\ldots,f_r)$
to be the cochain complex of graded modules
\begin{equation}  \label{RGam_f}
\RGam_{(f_1,\ldots,f_r)}(M) \, := \, \bigl[ \, M \xra{d^0} \prod_{1 \leq i_0 \leq r} M_{f_{i_0}}  \xra{d^1}  
\prod_{1 \leq i_0 < i_1 \leq r} M_{f_{i_0}f_{i_1}}  \xra{d^2} \ldots \xra{d^{r-1}} 
M_{f_1\ldots f_r} \, \bigr]
\end{equation}
whose differentials are defined by 
$d^m := \sum_{j=0}^m (-1)^j d^m_j$, where $d^m_j$ is the direct product of the canonical maps $d^m_j : A_{i_0\ldots \hat{i_j} \ldots i_m} \ra A_{i_0\ldots i_m}$.
Here, the first term $M$ is put in cohomological degree $0$. 

For a cochain complex $M \in \Ch(\Gr(A))$ of graded modules, we define $\RGam_{(f_1,\ldots,f_r)}(M)$ to be the total complex of the double complex $C^{p,q} = \RGam_{(f_1,\ldots,f_r)}(M^p)^q$.
\edf




The functor $M \mapsto \RGam_{(f_1,\ldots,f_r)}(M)$ on $\Ch(\Gr(A))$ is exact, and hence descend to a functor at the level of derived categories. Moreover, if we let $I := (f_1,\ldots,f_r)$ be the graded ideal generated by the elements $f_1$, then this functor has image inside the full subcategory $\cD_{\IIT}(\Gr (A))$. Thus, this gives a functor
\begin{equation}  \label{RGI_def}
\RGI := \RGam_{(f_1,\ldots,f_r)} \, : \, \cD(\Gr(A)) \raq \cD_{\IIT}(\Gr (A))
\end{equation}

Moreover, the map $\epsilon_M : \RGam_{(f_1,\ldots,f_r)}(M) \ra M$ defined by projecting to the first component of \eqref{RGam_f} gives rise to a natural transformation
\begin{equation}  \label{RGam_counit}
\epsilon \, : \, \iota \circ \RGI \, \Rightarrow \, \id
\end{equation}  
where $\iota : \cD_{\IIT}(\Gr (A)) \ra \cD(\Gr(A))$ is the inclusion functor. 
Then we have

\bthm   \label{RGam_right_adj}
The functor \eqref{RGI_def} is a right adjoint to the inclusion $\iota : \cD_{\IIT}(\Gr (A)) \ra \cD(\Gr(A))$, with counit given by \eqref{RGam_counit}. In particular the functor \eqref{RGI_def} depends only on the graded ideal $I$.
\ethm

The cone of $\epsilon_M$ is homotopic to the kernel of $\epsilon_M$ shifted by $1$, which is given by the following
\bdf  \label{Cech_cplx_def}
The \emph{\v{C}ech complex} of a graded module $M$ with respect to a tuple $(f_1,\ldots,f_r)$ of homogeneous elements
is the cochain complex of graded modules
\begin{equation}  \label{Ce_f}
\CI(M) \, = \, \Ce_{(f_1,\ldots,f_r)}(M) \, := \, \bigl[ \, \prod_{1 \leq i_0 \leq r} M_{f_{i_0}}  \xraq{-d^1}  
\prod_{1 \leq i_0 < i_1 \leq r} M_{f_{i_0}f_{i_1}}  \xraq{-d^2} \ldots \xraq{-d^{r-1}} 
M_{f_1\ldots f_r} \, \bigr]
\end{equation}
given as a subcomplex of \eqref{RGam_f}, 
shifted by one.
As in Definition \ref{RGam_f_def}, this definition can be extended to cochain complexes $M \in \Ch(\Gr(A))$ by taking the total complex.
\edf

Thus, for each $M \in \cD(\GrA)$, there is an exact triangle
\begin{equation}  \label{RGam_Ce_seq}
\ldots \raq \RGI(M) \xraq{\epsilon_M} M \xraq{\eta_M} \CI(M) \xraq{\delta_M} \RGI(M)[1] \raq \ldots 
\end{equation}
where $\eta_M = -d^0$, the negative of the first differential in \eqref{RGam_f},
and $\delta_M$ is the inclusion.

The exact triangle \eqref{RGam_Ce_seq} turns out to be the decomposition triangle associated to a semi-orthogonal decomposition. To this end, we recall the following
\blm  \label{I_triv_equiv}
For any $M \in \cD(\GrA)$, we have $\RGI(M) \simeq 0$ if and only if $\RHomcom_A(\RGI(A),M) \simeq 0$. 
\elm

\bdf
We say that $M \in \cD(\GrA)$ is \emph{$I$-trivial} if we have $\RGI(M) \simeq 0$, or equivalently $\RHomcom_A(\RGI(A),M) \simeq 0$ by Lemma \ref{I_triv_equiv}. Denote by $\cD_{\ITR}(\GrA) \subset \cD(\GrA)$ the full subcategory consisting of $I$-trivial objects.
\edf

\bpp
There is a semi-orthogonal decomposition
\begin{equation}  \label{local_cohom_SOD}
\cD(\GrA) \, = \, \langle \, \cD_{\ITR}(\GrA) \, , \,  \cD_{\IIT}(\GrA) \, \rangle 
\end{equation}
whose associated decomposition triangle is given by \eqref{RGam_Ce_seq}.
\epp

We are mostly interested in the special case $I = I^+ := A_{>0} \cdot A$, where the corresponding semi-orthogonal decomposition \eqref{local_cohom_SOD} is closely related to the projective space $X^+ = \Proj^+(A) := \Proj(A_{\geq 0})$.
Similarly, the semi-orthogonal decomposition \eqref{local_cohom_SOD} for $I = I^- := A_{<0} \cdot A$ is closely related to the projective space $X^- = \Proj^-(A) := \Proj(A_{\leq 0})$.

A first instance of this relation is the interpretation of $\Ce_{I^+}(M)$ in terms of a derived pushforward functor. Namely, consider the map $\pi^+ : \Proj^+(A) \ra \Spec(A_0) =: Y$, then for each $M \in \cD(\GrA)$ and for each $i \in \bZ$, there is a canonical isomorphism in $\cD(A_0)$:
\begin{equation}  \label{Ce_i_RGam}
\Ce_{I^+}(M)_i \, \cong \, {\bm R}\pi^+_*\,  \widetilde{M(i)} 
\end{equation}
where we denote by $\widetilde{M} \in \QCoh(X^+)$ the quasi-coherent sheaf on $X^+$ associated to $M \in \GrA$.
Indeed, \eqref{Ce_i_RGam} follows from the usual way of computing cohomology of quasi-coherent sheaves from the \v{C}ech complex (see \cite[(4.10)]{Yeu20c}).
This gives a proof of the following two results (see \cite[Lemma 4.12, 4.13]{Yeu20c} for details):
\blm  \label{RGam_Dbcoh}
If $A$ is Noetherian, then for any $M \in \Dbcoh(\GrA)$, we have $\Ce_{I^+}(M)_i \in \Dbcoh(A_0)$ and $\RGam_{I^+}(M)_i \in \Dbcoh(A_0)$ for each weight $i \in \bZ$.
\elm

\blm  \label{local_cohom_weight_bounded}
If $A$ is Noetherian, then for any $M \in \Dbcoh(\GrA)$, there exists $c^+, c^- \in \bZ$ such that 
$\RGam_{I^+}(M)_i \simeq 0$ for all $i \geq c^+$ and $\RGam_{I^-}(M)_i \simeq 0$ for all $i \leq c^-$.
\elm

Serre's equivalence gives another relation. We will use the version of Serre's equivalence proved in \cite{Yeu20c}, where the usual condition that $A$ be generated by $A_1$ over $A_0$ is replaced by the following condition:
\bdf  \label{frac_Cartier}
A $\bZ$-graded ring $A$ is said to be \emph{positively $\tfrac{1}{d}$-Cartier}, for an integer $d > 0$, if 
the canonical map $\widetilde{A(di)} \otimes_{\cO_{X^+}} \widetilde{A(dj)} \ra \widetilde{A(di+dj)}$ in $\QCoh(X^+)$ is an isomorphism for all $i,j \in \bZ$.
In the case $d = 1$, we simply say that $A$ is \emph{positively Cartier}.
\edf

For example, if $A_{\geq 0}$ is generated over $A_0$ by homogeneous elements $f_1,\ldots,f_p$ of positive degrees $d_i := \deg(f_i) > 0$, then it can be shown that $A$ is positively $\tfrac{1}{d}$-Cartier for any $d > 0$ that is divisible by each of $d_i$ (see \cite[Lemma 3.5]{Yeu20c}).

Denote by $\Tor^+ := \IpIT \subset \GrA$ the full subcategory consisting of $(I^+)^{\infty}$-torsion modules in the sense of Definition \ref{I_infty_torsion_def}, and denote by $\QpGr(A)$ the Serre quotient $\QpGr(A) := \GrA / \Torp$. Then we have (see \cite[Theorem 3.15]{Yeu20c})
\bthm  \label{Serre_equiv}
Suppose that $A$ is Noetherian and positively Cartier, then there is an equivalence of categories 
\begin{equation*}  
\begin{tikzcd}
\!\,^0 \! \cL^{+} \,:\, \QCoh(X^+) \ar[r, shift left]
& \QpGr(A) \, : \, (-)^{\sim}  \ar[l, shift left]
\end{tikzcd}
\end{equation*}
\ethm

\brm  \label{stacky_Proj}
When $A$ is not positively Cartier, the category $\QpGr(A)$ also admit a description in terms of quasi-coherent sheaves on a stacky projective space. 
Namely, one can show that the map of $\bG_m$-equivariant schemes $\Spec(A) \ra \Spec(A_{\geq 0})$ induces an isomorphism on the $\bG_m$-invariant open subschemes
\begin{equation*}
\Spec(A) \setminus \Spec(A/I^+) \xraq{\cong} \Spec(A_{\geq 0}) \setminus \Spec(A_0)
\end{equation*}

If we denote this $\bG_m$-equivariant scheme by $W^{ss}(+)$, and let $\, \PProj^+(A)$ be the quotient stack $[W^{ss}(+)/\bG_m]$, then we have an equivalence
$
\QpGr(A)  \simeq  \QCoh(\PProj^+(A)) 
$.
See, e.g., \cite[Proposition 2.3]{AKO08} or 
\cite[Example 2.15]{HL15}.
\erm

Given a Serre quotient $\phi^* : \cC \ra \cC/\cT$, there are some general criteria laid out in \cite[Appendix A]{Yeu20c} which guarantees that the derievd category of the quotient $\cD(\cC/\cT)$ coincides with the Verdier quotient $\cD(\cC)/\cD_{\cT}(\cC)$ of the derived categories. 
If $A$ is Noetherian and positively Cartier, then Theorem \ref{Serre_equiv} asserts that $\QCoh(X^+)$ is a Serre quotient. The corresponding criteria can be easily verified, and we have

\bpp[\cite{Yeu20c}, Proposition 3.20]
The functor $\phi^* : \cD(\GrA) \ra \cD(\QpGr(A))$ has a fully faithful right adjoint ${\bm R}\phi_* : \cD(\QpGr(A)) \ra \cD(\GrA)$, which induces a semi-orthogonal decomposition
\begin{equation*}
\cD(\GrA) \, = \, \langle \,  {\bm R}\phi_* ( \cD(\QpGr(A))) \,   , \, \cD_{\Torp}(\GrA) \, \rangle
\end{equation*}
Comparing with \eqref{local_cohom_SOD}, we see that the following is a pair of inverse exact equivalences:
\begin{equation}  \label{QpGrA_IpTR_equiv}
\begin{tikzcd}
\phi^* \, : \, \Dsuit_{\IpTR}(\GrA) \ar[r, shift left]
&  \Dsuit(\QpGr(A)) \, : \, {\bm R}\phi_* \ar[l, shift left]
\end{tikzcd}
\end{equation}
\epp

The exact equivalence \eqref{QpGrA_IpTR_equiv} restricts to certain bounded coherent subcategories. One has to be careful that the relevant subcategory of $\cD_{\IpTR}(\GrA)$ is not given by $\cD_{\IpTR}(\GrA) \cap \Dbcoh(\GrA)$. Instead, one considers the following
\bdf  \label{IpTR_coh_def}
Denote by  $\gr(A) \subset \GrA$ the full subcategory of finitely generated graded modules.
Let $\qpgr(A) \subset \QpGr(A)$ be the essentially image of $\gr(A)$ under $\phi^* : \GrA \ra \QpGr(A)$.
For each $\spadesuit \in \{\, \, \, ,+,-,b\}$, 
\begin{enumerate}
	\item let $\Dsuit_{\coh}(\QpGr(A)) \subset \Dsuit(\QpGr(A))$ be the full subcategory consisting of complexes whose cohomology lies in $\qpgr(A)$; and
	\item let $\Dsuit_{{\rm coh}(\IpTR)}(\GrA) \subset \Dsuit_{\IpTR}(\GrA)$ be the essential image of $\Dsuit_{{\rm coh}}(\GrA)$ under the functor $\Ce_{I^+} : \Dsuit(\GrA) \ra \Dsuit_{\IpTR}(GrA)$.
\end{enumerate}
\edf

\bpp [\cite{Yeu20c}, Corollary 3.26] \label{coh_IpTR_QpGr}
For $\spadesuit \in \{-,b\}$, the equivalence \eqref{QpGrA_IpTR_equiv} restricts to give an exact equivalence 
\begin{equation*} 
\begin{tikzcd}
\phi^* \, : \, \Dsuit_{\coh(\IpTR)}(\GrA) \ar[r, shift left]
&  \Dsuit_{\coh}(\QpGr(A)) \, : \, {\bm R}\phi_* \ar[l, shift left]
\end{tikzcd}
\end{equation*}
\epp

Thus, if $A$ is positively Cartier, then the subcategory $\cD^b_{\coh(\IpTR)}(\GrA)$ is equivalent to the usual bounded derived category of coherent sheaves on $X^+ = \Proj^+(A)$. In general, it is equivalent to that of the stacky projective space $\PProj^+(A)$ by Remark \ref{stacky_Proj}.

\section{Weight truncation}

In this section, we use freely the language of modules over small pre-additive categories. The reader is referred to Appendix \ref{app_mod_preadd} for details.

Given any $\bZ$-graded ring $A$,
let $\cF = \cF_A$ be the pre-additive category with object set $\Ob(\cF) = \bZ$, and Hom spaces $\cF(i,j) := A_{i-j}$. Compositions in $\cF$ are defined by multiplication in $A$ in the obvious way.
A right $\cF$-module is nothing but a graded (right) $A$-module. More precisely, there is an equivalence of abelian categories
\begin{equation}  \label{Fmod_GrA_1}
(-)^{\sharp} \, : \, \GrA \xraq{\simeq} \Mod(\cF) \, , \qquad \quad (M^{\sharp})_i := M_i
\end{equation}
whose inverse will be denoted as $(-)^{\flat} : \Mod(\cF) \xra{\simeq} \GrA$.

Since $A$ is assumed to be commutative, any graded module $M \in \GrA$ in fact induces an $\cF$-bimodule. In other words, there is an additive functor
\begin{equation}  \label{tilde_F_bimod}
\Gr(A) \raq \Mod(\cF^e) \, , \quad M \mapsto \widetilde{M} \, , \quad \text{where} \quad \,_j\widetilde{M}_i := M_{i-j}
\end{equation}
which recovers the functor \eqref{Fmod_GrA_1} by restricting to $\,_0\widetilde{M}_{*}$. 

The graded tensor products and graded Hom spaces between graded modules can be expressed naturally in terms of $\cF$-bimodules. Namely, for any $M,N \in \GrA$, there are natural isomorphisms of $\cF$-bimodules 
\begin{equation}  \label{cF_A_tensor_Hom}
\widetilde{M} \otimes_{\cF} \widetilde{N} \, \cong \, \widetilde{M \otimes_A N} 
\qquad \text{ and } \qquad 
\Hom_{\cF}(\widetilde{M} , \widetilde{N}) \, \cong \, \widetilde{\Homcom_A(M,N)}
\end{equation}
In particular, if we only remember the right $\cF$-module structure of $\widetilde{M}$, then we  have
\begin{equation}  \label{sharp_tensor_tilde}
M^{\sharp} \otimes_{\cF} \widetilde{N} \, \cong \, (M \otimes_A N)^{\sharp} 
\end{equation}

Since $A$ is assumed to be commutative, the pre-additive category $\cF = \cF_A$ admits an involution, \ie, it comes equipped with an isomorphism $\cF \cong \cF^{\op}$ of pre-additive categories, defined by $i \mapsto -i$ on objects, and $\cF(i,j) = A_{i-j} = \cF^{\op}(-i,-j)$ on Hom sets. In fact, the commutativity of $A$ is equivalent to the fact that this assignment $\cF \ra \cF^{\op}$ is a functor.
This involution induces an isomorphism of categories
\begin{equation*}
(-)^{\tau} \, : \, \Mod(\cF) \xraq{\cong} \Mod(\cF^{\op}) \, , \qquad \,_i(M^{\tau}) := M_{-i}
\end{equation*}
whose (strict) inverse will also be denoted as $(-)^{\tau}$. 

For any integer $a \in \bZ$, let $\Fga$ be the full subcategory of $\cF$ on the subset $\Ob(\Fga) = \bZ_{\geq a} \subset \bZ = \Ob(\cF)$. 
Define $\Fla \subset \cF$ in the similar way. 
We will also write $\cF_{\leq \infty} = \cF = \cF_{\geq -\infty}$.
For any $-\infty \leq a' \leq a$, denote by $(-)_{\geq a} : \Mod(\cF_{\geq a'}) \ra \Mod(\cF_{\geq a})$ the restriction functor.

Notice that the involution on $\cF$ restricts to an isomorphism $(\Fga)^{\op} \cong \Flma$. 
As a result, there is again an isomorphism of categories
\begin{equation}  \label{transpose_Fga}
(-)^{\tau} \, : \, \Mod(\Fga) \xraq{\cong} \Mod((\Flma)^{\op}) \, , \qquad \,_i(M^{\tau}) := M_{-i}
\end{equation}
whose inverse will also be denoted as $(-)^{\tau}$. Similarly, the functor $(M^{\tau})_i := \!\,_{-i}M$ gives an isomorphism of categories $(-)^{\tau} : \Mod((\Fga)^{\op}) \xra{\cong} \Mod(\Flma)$ whose inverse will also be denoted as $(-)^{\tau}$.

From now on, we fix an integer $w \in \bZ$.
Notice that the inclusion functor  $\Fgw \rinto \cF$ induces a three-way adjunction
\begin{equation} \label{Fga_F_adj}
\begin{tikzcd}
\Mod(\Fgw) \ar[rr, bend left, "-\otimes_{\Fgw} \cF"] \ar[rr, bend right, "\Hom_{\Fgw}(\cF\text{,}-)"']
& &  \Mod(\cF) \ar[ll, "(-)_{\geq w}"']
\end{tikzcd}
\end{equation} 

In fact, under the equivalence $\Mod(\cF) \simeq \GrA$, the right-pointing functor $- \otimes_{\Fgw} \cF$ on the top may be characterized as the unqiue cocontinuous functor satisfying
\begin{equation}  \label{weight_tensor_free}
\parbox{40em}{For each $i \geq w$, the functor $(- \otimes_{\Fgw} \cF)^{\flat} : \Mod(\Fgw) \ra \Gr(A)$ sends the free module $\!\,_{i}\Fgw$ to the free graded module $A(-i)$}
\end{equation}

Since the inclusion functor $\Fgw \ra \cF$ is fully faithful, we have
\begin{equation}  \label{tensor_restrict_Fgw}
\begin{split}
(M \otimes_{\Fgw} \cF)_{\geq w} \, &\cong \, M \otimes_{\Fgw} \Fgw \, \cong \, M \\
\Hom_{\Fgw}(\cF, M)_{\geq w} \, &\cong \, \Hom_{\Fgw}(\Fgw , M) \, \cong \, M
\end{split}
\end{equation}

One can also relate $\Fgw$-modules to graded $A_{\geq 0}$-modules. Indeed, for any $\cM \in \Mod(\Fgw)$, the assignment $M_i := \cM_i$ defines a graded $A_{\geq 0}$-module, which will be denoted as $\cM|_{A_{\geq 0}}$. This can be used to characterize finite generated $\Fgw$-modules:

\blm  \label{Fga_fg_module}
Suppose that $A$ is Noetherian, then for any $\cM \in \Mod(\Fgw)$, the followings are equivalent:
\begin{enumerate}
	\item $\cM$ is a finitely generated $\Fgw$-module in the sense of Definition \ref{fin_gen_mod};
	\item there exists a finitely generated graded $A$-module $M$ such that $\cM \cong M^{\sharp}_{\geq w} := (M^{\sharp})_{\geq w}$;
	\item $\cM|_{A_{\geq 0}}$ is a finitely generaeted graded $A_{\geq 0}$-module.
\end{enumerate}
\elm

\bpf
For $(1) \Rightarrow (2)$, simply notice that if there is an epimorphism $\oplus_{j = 1}^m \, \,_{i_j}\Fgw \ronto \cM$, then applying $-\otimes_{\Fgw} \cF$, we have by \eqref{weight_tensor_free} an epimorphism 
$\oplus_{j = 1}^m \, A(-i_j) \ronto (\cM \otimes_{\Fgw} \cF)^{\flat}$, which shows that $M := (\cM \otimes_{\Fgw} \cF)^{\flat}$ is a finitely generated graded $A$-module. Moreover, it satisfies $\cM \cong M^{\sharp}_{[\geq w]}$ by \eqref{tensor_restrict_Fgw}.

For $(2) \Rightarrow (3)$, we use the Noetherian condition. By Proposition \ref{Noeth_gr_ring}, $A_0$ is Noetherian, and each $A_i$ is finitely generated as a module over $A_0$. 
Moreover, there exists integers $d>0$ and $N_0 > 0$ such that whenever $N \geq N_0$, we have $A_{N} = A_{d} \cdot A_{N-d}$ (see, for example, \cite[Lemma 3.2]{Yeu20c}). Now if $M$ is generated by $\xi_1,\ldots,\xi_r$ over $A$, then for any $j \geq N_0 + \max_i\{ \deg(\xi_i) \}$, we have $M_j = A_d \cdot M_{j-d}$. 
Since each $M_i$ is clearly finitely generated over $A_0$, the same is true for $\oplus_{w \leq j \leq N_0 + \max_i\{ \deg(\xi_i) \}} M_i$, which then generates $M^{\#}_{[\geq w]}|_{A_{\geq 0}}$ as an $A_{\geq 0}$-module.

The implication $(3) \Rightarrow (1)$ follows directly from the definitions.
\epf

\bcor  \label{Fgw_Noeth}
If the $\bZ$-graded ring $A$ is Noetherian, then the small pre-additive category $\Fgw$ is Noetherian.
\ecor

\bpf
By the characterization of finitely generated right $\Fgw$-modules in Lemma \ref{Fga_fg_module}(3), we see that $\Fgw$ is right Noetherian. The same argument also shows that $\cF_{\leq -w}$ is right Noetherian. By the isomorphism of categories $(-)^{\tau} : \Mod((\Fgw)^{\op}) \xra{\cong} \Mod(\cF_{\leq -w})$, we see that $\Fgw$ is left Noetherian as well.
\epf

Now we take the derived functors of the functors appearing in \eqref{Fga_F_adj}. Under the identification $\cD(\cF) \simeq \cD(\GrA)$, we denote these derived functors as
\begin{equation*}
\begin{split}
\mathscr{L}_{[\geq w]} \, &: \, \cD(\Fgw) \raq \cD(\GrA) \, , \qquad \mathscr{L}_{[\geq w]}(\cM) \, := \, (\cM \otimes^{{\bm L}}_{\Fgw} \cF)^{\flat} \\
\mathscr{R}_{\{\geq w\}} \, &: \, \cD(\Fgw) \raq \cD(\GrA) \, , \qquad \mathscr{R}_{\{\geq w\}}(\cM) \, := \,
(\RHom_{\Fgw}(\cF , \cM) ) ^{\flat}
\end{split}
\end{equation*}

\bdf
Let $\cD_{<w}(\GrA) \subset \cD(\GrA)$ be the full subcategory consisting of $M \in \cD(\GrA)$ such that $M_i \simeq 0$ for all $i \geq w$.
\edf

Then we have a recollement
\begin{equation}  \label{weight_recollement}
\begin{tikzcd}
\cD_{<w}(\GrA) \ar[rr, "\iota"]
&& \cD(\GrA)  \ar[rr, "(-)^{\sharp}_{\geq w}"]  \ar[ll, bend right, "\cL_{<w}"'] \ar[ll, bend left, "\cR_{<w}"']
&& \cD(\cF_{\geq w})  \ar[ll, bend right, "\mathscr{L}_{[\geq w]}"'] \ar[ll, bend left, "\mathscr{R}_{\{\geq w\}}"']
\end{tikzcd}
\end{equation}
where we have written $M^{\sharp}_{\geq w} := (M^{\sharp})_{\geq w}$.

Indeed, if we denote by $\cL_{[\geq w]}$ and $\cR_{\{\geq w\}}$ the endofunctors on $\cD(\GrA)$ given by the compositions
\begin{equation*}
\cL_{[\geq w]}(M) \, := \, \mathscr{L}_{[\geq w]}(M^{\#}_{\geq w})  \qquad \text{and} \qquad 
\cR_{\{\geq w\}}(M) \, := \, \mathscr{R}_{\{\geq w\}}(M^{\#}_{\geq w})
\end{equation*}
then the functors $\cL_{< w}$ and $\cR_{<w}$ are defined by the exact triangles
\begin{equation}  \label{weight_trunc_exact_tri}
\begin{split}
\ldots \raq \cL_{[\geq w]}(M) \raq &M \raq \cL_{< w}(M) \raq  \cL_{[\geq w]}(M)[1] \raq \ldots \\
\ldots \raq \cR_{< w}(M) \raq &M \raq \cR_{\{\geq w\}}(M) \raq  \cR_{<w}(M)[1] \raq \ldots
\end{split}
\end{equation}

We are mostly interested in the functors $\mathscr{L}_{[\geq w]}$, $\cL_{[\geq w]}$ and $\cL_{<w}$ instead of their counterpart $\mathscr{R}_{\{\geq w\}}$, $\cR_{\{\geq w\}}$ and $\cR_{<w}$. In particular, we emphasize the following

\bdf
Let $\cD_{[\geq w]}(\GrA) \subset \cD(\GrA)$ be the essential image of the fully faithful functor $\mathscr{L}_{[\geq w]}  : \cD(\Fgw) \ra \cD(\GrA)$.
\edf

Alternatively, the subcategory $\cD_{[\geq w]}(\GrA) \subset \cD(\GrA)$ may be characterized as follows:
\blm  \label{D_geq_w_gen}
$\cD_{[\geq w]}(\GrA)$ is the smallest strictly full triangulated subcategory containing the objects $A(-i)$ for $i \geq w$ and is closed under small coproducts.
Therefore, we have
\begin{equation*}
\cD_{[\geq w]}(\GrA) \cap \cD(\GrA)_c \, = \, \cD_{[\geq w]}(\GrA)_c  \, = \, {\rm EssIm}(\,\scL_{[\geq w]} : \cD_{\perf}(\Fgw) \ra \cD(\GrA) \,)
\end{equation*}
where the subscript $(-)_c$ denotes the full subcategory of compact objects.
\elm

\bpf
The first statement follows from \eqref{weight_tensor_free} and the fact that $\scL_{[\geq w]}$ preserves small coproducts. For the second statement, the first equality is standard (see, e.g., \cite[Lemma 2.2]{Nee92} or \cite[Theorem 5.3]{Rou08}). The second equality follows from the standard fact (see, e.g., \eqref{cpt_equal_perf}) that $\cD_{\perf}(\Fgw) = \cD(\Fgw)_c$.
\epf

It follows immediately from the recollement \eqref{weight_recollement} that there is a semi-orthogonal decomposition
\begin{equation}  \label{weight_SOD}
\cD(\GrA) = \langle \, \cD_{< w}(\GrA) \, , \, \cD_{[\geq w]}(\GrA) \, \rangle
\end{equation}

\brm
The notation $\cD_{[\geq w]}$ is meant to convey the idea that these are objects ``generated in weight $\geq w$"; while the notation $\cD_{< w}$ means that these are objects ``concentrated in weight $<w$".
\erm

\section{Local cohomology and weight truncation}

Now we study local cohomology under weight truncation. For any finitely generated graded ideal $I \subset A$, the objects $\RGI(A)$ and $\CI(A)$ in $\cD(\GrA)$ give rise to the objects $\!\,_{\geq w}\widetilde{\RGI(A)}_{\geq w}$ and $\!\,_{\geq w}\widetilde{\CI(A)}_{\geq w}$ in the derived category $\cD((\Fgw)^e)$ of $\Fgw$-bimodules.
Tensoring over these give rise to functors
\begin{equation}  \label{RGam_Ce_weight}
\begin{split}
\RGam_{I,\geq w}  &: \cD(\Fgw) \ra \cD(\Fgw) , \quad \cM  \mapsto  \cM \otimes^{{\bm L}}_{\Fgw} \widetilde{\RGI(A)}_{\geq w}  \cong 
((\mathscr{L}_{[\geq w]}\cM) \otimes_A^{{\bm L}} \RGI(A) )^{\sharp}_{\geq w} \\
\Ce_{I,\geq w}  &: \cD(\Fgw) \ra \cD(\Fgw) , \quad \cM  \mapsto  \cM \otimes^{{\bm L}}_{\Fgw} \widetilde{\CI(A)}_{\geq w}  \cong 
((\mathscr{L}_{[\geq w]}\cM) \otimes_A^{{\bm L}} \CI(A) )^{\sharp}_{\geq w}
\end{split}
\end{equation}
where the last isomorphisms on each line is obtained by applying \eqref{tensor_two_step} and \eqref{sharp_tensor_tilde}.

We are mostly interested in the case $I = I^+$, where these functors behave very similarly to the corresponding ones on $\cD(\GrA)$ (see Proposition \ref{RGam_Ce_weight_descend} and Theorem \ref{RGam_Ce_weight_ortho} below). In fact, these properties are formal consequences of the following easy
\blm  \label{weight_subset_tor}
We have $\cD_{<w}(\GrA) \subset \cD_{\Torp}(\GrA)$.
\elm

which can be used to prove the following two results:

\blm  \label{Torp_weight_lem}
If $M \in \cD_{\Torp}(\GrA)$ then $\cL_{[\geq w]}(M) \in \cD_{\Torp}(\GrA)$.

If $\RGam_{I^+}(M) \in \cD_{<w}(\GrA)$ then $\cL_{[\geq w]}(M) \otimes_A^{{\bm L}} \RGam_{I^+}(A) \in \cD_{< w}(\GrA)$.
\elm

\bpf
For any $M \in \cD(\GrA)$, we have $\cL_{<w}(M) \in \cD_{\Torp}(\GrA)$ by Lemma \ref{weight_subset_tor}, so that the first statement follows directly from the exact triangle in the first row of \eqref{weight_trunc_exact_tri}.
For the second statement, apply $\RGam_{I^+}(-)$ to the same exact triangle, and 
notice that  $\cL_{<w}(M) \otimes_A^{{\bm L}} \RGam_{I^+}(A) \simeq \cL_{<w}(M) \in \cD_{<w}(\GrA)$ because $\cL_{<w}(M) \in \cD_{\Torp}(\GrA)$.
\epf

\bpp  \label{RGam_Ce_weight_descend}
The following two functors commute up to isomorphism of functors:
\begin{equation*}
\begin{tikzcd}
\cD(\GrA) \ar[r, "\RGam_{I^+}"] \ar[d, "(-)^{\sharp}_{\geq w}"'] & \cD(\GrA) \ar[d, "(-)^{\sharp}_{\geq w}"]
& &  \cD(\GrA) \ar[r, "\Ce_{I^+}"] \ar[d, "(-)^{\sharp}_{\geq w}"'] & \cD(\GrA) \ar[d, "(-)^{\sharp}_{\geq w}"] \\
\cD(\Fgw) \ar[r, "\RGam_{I^+,\geq w}"] & \cD(\Fgw) 
& &  \cD(\Fgw) \ar[r, "\Ce_{I^+,\geq w}"] & \cD(\Fgw)
\end{tikzcd}
\end{equation*}
\epp

\bpf
Given any $M \in \cD(\GrA)$, take the exact triangle in the first row of \eqref{weight_trunc_exact_tri}. By Lemma \ref{weight_subset_tor}, we have $\cL_{<w}(M) \in \cD_{\Torp}(\GrA)$. Applying $\Ce_{I^+}$ to this exact triangle, we have $\Ce_{I^+}(\cL_{[\geq w]} M) \cong \Ce_{I^+}(M)$ in $\cD(\GrA)$. Applying $(-)^{\sharp}_{\geq w}$ to this isomorphism gives the commutativity of the second square.
Similarly, applying $\RGam_{I^+}$ to the same exact triangle, we have an exact triangle 
\begin{equation*}
\ldots \raq \RGam_{I^+}(\cL_{[\geq w]}(M)) \raq \RGam_{I^+}(M) \raq \cL_{< w}(M) \raq  \RGam_{I^+}(\cL_{[\geq w]}(M))[1] \raq \ldots
\end{equation*}
Applying $(-)^{\#}_{\geq w}$ therefore gives an isomorphism $(\RGam_{I^+}(\cL_{[\geq w]}(M)))^{\sharp}_{\geq w} \cong (\RGam_{I^+}(M))^{\sharp}_{\geq w}$, proving
the commutativity of the first square.
\epf

%

\bdf
Let $\cD_{\IpTR}(\Fgw)$ and $\cD_{\Torp}(\Fgw)$ be the full subcategories of $\cD(\Fgw)$ defined by
\begin{equation*}
\begin{split}
\cD_{\IpTR}(\Fgw) \,&:= \, \{ \, \cM \in \cD(\Fgw) \, | \, \RGam_{I^+,\geq w}(\cM) \simeq 0 \, \} \\
\cD_{\Torp}(\Fgw) \,&:= \, \{ \, \cM \in \cD(\Fgw) \, | \, \Ce_{I^+,\geq w}(\cM) \simeq 0 \, \}
\end{split}
\end{equation*}
\edf

%

Then we have the following

\bthm  \label{RGam_Ce_weight_ortho}
For $I = I^+$, the functors \eqref{RGam_Ce_weight} form a semi-orthogonal pair of idempotents, in the sense that the followings hold:
\begin{enumerate}
	\item for any $\cM \in \cD(\Fgw)$, we have $\RGam_{I^+,\geq w}(\cM) \in \cD_{\Torp}(\Fgw)$ 
	and $\Ce_{I^+,\geq w}(\cM) \in \cD_{\IpTR}(\Fgw)$;
	\item if $\cM \in \cD_{\Torp}(\Fgw)$ then $\RGam_{I^+,\geq w}(\cM) \cong \cM$; 
	\item if $\cM \in \cD_{\IpTR}(\Fgw)$ then $\Ce_{I^+,\geq w}(\cM) \cong \cM$;
	\item there is a semi-orthogonal decomposition $\cD(\Fgw) = \langle \, \cD_{\IpTR}(\Fgw) \, , \, \cD_{\Torp}(\Fgw) \, \rangle$.
\end{enumerate}
\ethm

\bpf
Statement (1) follows immediately from Proposition \ref{RGam_Ce_weight_descend}. For statements (2) and (3), notice that every $\cM \in \cD(\Fgw)$ can be written as $\cM \cong M^{\sharp}_{\geq w}$ for some $M \in \cD(\GrA)$. For example one can take $M = \scL_{[\geq w]}(\cM)$.
For this $M$, take the exact triangle \eqref{RGam_Ce_seq} for $I = I^+$, and apply $(-)^{\sharp}_{\geq w}$ to the triangle.
In view of Proposition \ref{RGam_Ce_weight_descend}, we have an exact triangle
\begin{equation*}
\ldots \raq \RGam_{I^+,\geq w}(\cM) \xraq{\epsilon_{\cM}} \cM \xraq{\eta_{\cM}} \Ce_{I^+,\geq w}(\cM) \xraq{\delta_{\cM}} \RGam_{I^+,\geq w}(\cM)[1] \raq  \ldots 
\end{equation*}
which immediately shows (2) and (3). In fact, because of (1), this exact triangle also establishes the decomposition for (4), so that it suffices to show the semi-orthogonality 
$\cD_{\Torp}(\Fgw) \perp \cD_{\IpTR}(\Fgw)$. 
Thus, let $\cM \in \cD_{\Torp}(\Fgw)$ and $\cN \in \cD_{\IpTR}(\Fgw)$, then by (2) and (3), we may write $\cM \cong M^{\sharp}_{\geq w}$ and $\cN \cong N^{\sharp}_{\geq w}$ for some $M \in \cD_{\Torp}(\GrA)$ and $N \in \cD_{\IpTR}(\GrA)$. By the adjunction \eqref{weight_recollement}, we have
\begin{equation*}
\Hom_{\cD(\Fgw)}( \cM , \cN ) \, \cong \, \Hom_{\cD(\GrA)} ( \cL_{[\geq w]}(M), N)
\end{equation*}
which is zero because of the first statement of Lemma \ref{Torp_weight_lem}. 
\epf

Combined with \eqref{weight_SOD}, this gives the first statement of the following

\bthm  \label{three_term_SOD}
There is a semi-orthogonal decomposition 
\begin{equation*}
\cD(\GrA) \, = \, \langle \, \cD_{<w}(\GrA) \, , \, \mathscr{L}_{[\geq w]} ( \cD_{\IpTR}(\Fgw) )  \, , \,  \mathscr{L}_{[\geq w]} ( \cD_{\Torp}(\Fgw) ) \, \rangle
\end{equation*}
where the functor $\mathscr{L}_{[\geq w]} : \cD(\Fgw) \ra \cD(\GrA)$ is fully faithful. 
Moreover, the latter two semi-orthogonal components can be identified as
\begin{equation}  \label{SOD_comp_descr_1}
\begin{split}
\mathscr{L}_{[\geq w]} ( \cD_{\IpTR}(\Fgw) ) \, &= \, \{ \, M \in \cD_{[\geq w]}(\GrA) \, | \, \RGam_{I^+}(M) \in \cD_{<w}(\GrA) \, \} \\
\mathscr{L}_{[\geq w]} ( \cD_{ \Torp}(\Fgw) ) \, &= \,
\cD_{\Torp,[\geq w]}(\GrA) \, := \,  \cD_{[\geq w]}(\GrA) \cap \cD_{\Torp}(\GrA)
\end{split}
\end{equation}
\ethm

\bpf
Only the identifications \eqref{SOD_comp_descr_1} of the semi-orthogonal components needs proof. 
It is clear that all the subcategories in \eqref{SOD_comp_descr_1} lie in $\cD_{[\geq w]}(\GrA)$. Moreover, given $M \in \cD_{[\geq w]}(\GrA)$, then by Proposition \ref{RGam_Ce_weight_descend}, we have $M^{\sharp}_{\geq w} \in \cD_{\IpTR}(\Fgw)$ if and only if $(\RGam_{I^+}(M))^{\sharp}_{\geq w} = 0$. The latter condition is precisely $\RGam_{I^+}(M) \in \cD_{<w}(\GrA)$, hence proving the first row of \eqref{SOD_comp_descr_1}.
For the second row, we similarly observe that, given any $M \in \cD_{[\geq w]}(\GrA)$, then by Proposition \ref{RGam_Ce_weight_descend}, we have $M^{\sharp}_{\geq w} \in \cD_{\Torp}(\Fgw)$ if and only if $\Ce_{I^+}(M) \in \cD_{<w}(\GrA)$. But $\Ce_{I^+}(M)$ is always in $\cD_{\IpTR}(\GrA)$, so that the latter is true if and only if $\Ce_{I^+}(M) = 0$.
\epf

Unravelling the definitions, we see that this semi-orthogonal decomposition decomposes every $M \in \cD(\GrA)$ into the diagram
\begin{equation}  \label{triple_SOD_terms}
\begin{tikzcd} [row sep = 12, column sep = 15]
\cL_{[\geq w]} \RGam_{I^+}(M)  \ar[rr, "\cL_{[\geq w]} (\epsilon_M)"] 
& &  \cL_{[\geq w]} M   \ar[rr, "\text{counit}"] 
\ar[ld, " \cL_{[\geq w]}(\eta_M)"]
& & M \ar[ld, "\text{unit}"] \\
& \cL_{[\geq w]} \Ce_{I^+}(M)  \ar[ul, "\text{[1]}" description]
& & \cL_{<w}(M) \ar[ul, "\text{[1]}" description]
\end{tikzcd}
\end{equation}
with the decomposition terms 
\begin{enumerate}
	\item $\cL_{<w}(M) \in  \cD_{<w}(\GrA)$,
	\item $\cL_{[\geq w]} \Ce_{I^+}(M) \in  \mathscr{L}_{[\geq w]} ( \cD_{\IpTR}(\Fgw) )$, and 
	\item $\cL_{[\geq w]} \RGam_{I^+}(M) \in \mathscr{L}_{[\geq w]} ( \cD_{\Torp}(\Fgw) ) = \cD_{\Torp,[\geq w]}(\GrA)$
\end{enumerate}

We also have the following two characterizations of semi-orthogonal components
\blm  \label{Torp_SOD}
There is a semi-orthogonal decomposition 
\begin{equation*}
\cD_{\Torp}(\GrA) \, = \, \langle \, \cD_{<w}(\GrA) \, , \, \cD_{\Torp,[\geq w]}(\GrA) ) \, \rangle
\end{equation*}
\elm

\bpf
For any $M \in \cD_{\Torp}(\GrA)$, take the exact sequence in the first row of \eqref{weight_trunc_exact_tri}. Since $\cL_{<w}M \in \cD_{<w}(\GrA) \subset \cD_{\Torp}(\GrA)$, we have $\cL_{[\geq w]}(M) \in \cD_{\Torp}(\GrA)$ as well. This proves the claimed decomposition. Orthogonality was already shown in Theorem \ref{three_term_SOD}.
\epf

\blm  \label{first_two_terms_SOD}
The composite of the first two components in Theorem \ref{three_term_SOD} can be identified as
\begin{equation*}
\langle \, \cD_{<w}(\GrA) \, , \, \mathscr{L}_{[\geq w]} ( \cD_{\IpTR}(\Fgw) ) \, \rangle \,
= \, \{ \, M \in \cD(\GrA) \, | \, \RGam_{I^+}(M) \in \cD_{<w}(\GrA) \, \}  
\end{equation*}
\elm

\bpf
An object $M \in \cD(\GrA)$ lies in $\langle \, \cD_{<w}(\GrA) \, , \, \mathscr{L}_{[\geq w]} ( \cD_{\IpTR}(\Fgw) ) \, \rangle$ if and only if the component $\cL_{[\geq w]} \RGam_{I^+}(M)$ in \eqref{triple_SOD_terms} vanishes. This is true if and only if $\RGam_{I^+}(M) \in \cD_{<w}(\GrA)$.
\epf

Our next goal is to show that the triangulated category $\cD_{\IpTR}(\Fgw)$ is in fact equivalent to $\cD_{\IpTR}(\GrA)$ (see Theorem \ref{I_triv_weight_equiv} below). First, notice that, by the second statement of Lemma \ref{Torp_weight_lem}, the functor $(-)^{\sharp}_{\geq w} : \cD(\GrA) \ra \cD(\Fgw)$ sends the subcategory $\cD_{\IpTR}(\GrA) \subset \cD(\GrA)$ to the subcategory $\cD_{\IpTR}(\Fgw) \subset \cD(\Fgw)$, so that we have an exact functor
\begin{equation}  \label{I_triv_weight_functor_1}
(-)^{\sharp}_{\geq w} \, : \, \cD_{\IpTR}(\GrA) \raq \cD_{\IpTR}(\Fgw)
\end{equation}

In the other direction, there is the exact functor
\begin{equation}  \label{I_triv_weight_functor_2}
\Ce_{I^+} \circ \mathscr{L}_{[\geq w]} \, : \, \cD_{\IpTR}(\Fgw) \raq \cD_{\IpTR}(\GrA)
\end{equation}

\bthm  \label{I_triv_weight_equiv}
The functors \eqref{I_triv_weight_functor_1} and \eqref{I_triv_weight_functor_2} are inverse equivalences.
\ethm

\bpf
First, notice that the functors \eqref{I_triv_weight_functor_1} and \eqref{I_triv_weight_functor_2} are restrictions of the composite adjunctions
\begin{equation*}
\begin{tikzcd}
\cD(\Fgw) \ar[r, shift left, "\mathscr{L}_{[\geq w]}"]
& \cD(\GrA) \ar[r, shift left, "\Ce_{I^+}"]  \ar[l, shift left, "(-)^{\sharp}_{\geq w}"]
& \cD_{\IpTR}(\GrA)  \ar[l, shift left, "\iota"]
\end{tikzcd}
\end{equation*}
and are therefore adjoints to each other.

The fact that the adjunction unit $\id \Rightarrow (-)^{\sharp}_{\geq w} \circ \Ce_{I^+} \circ \mathscr{L}_{[\geq w]}$  is an isomorphism on $\cD_{\IpTR}(\Fgw)$ is precisely statement (3) of Theorem \ref{RGam_Ce_weight_ortho}. 
This shows that, for any $M \in \cD(\GrA)$, the adjunction counit $\Ce_{I^+} ( \mathscr{L}_{[\geq w]} (M^{\sharp}_{\geq w})) \ra M$ becomes an isomorphism after applying $(-)^{\sharp}_{\geq w}$. In other words, its cone lies in $\cD_{<w}(\GrA) \subset \cD_{\Torp}(\GrA)$. If $M \in \cD_{\IpTR}(\GrA)$, then this cone also lies in $\cD_{\IpTR}(\GrA)$, which means it must be zero.
\epf

The semi-orthogonal decomposition in Theorem \ref{three_term_SOD} can be rewritten in the form
\begin{equation}  \label{SOD_no_F}
\cD(\GrA) \, = \, \langle \, \cD_{<w}(\GrA) \, , \, \cL_{[\geq w]} (\cD_{\IpTR}(\GrA)) \, , \,  \cD_{\Torp,[\geq w]}(\GrA)  \, \rangle
\end{equation}
where $\cL_{[\geq w]} (\cD_{\IpTR}(\GrA))$ is the essential image of the functor $\cL_{[\geq w]} : \cD_{\IpTR}(\GrA) \ra \cD(\GrA)$, which is fully faithful by Theorem \ref{I_triv_weight_equiv}.

\brm
Our discussion so far about local cohomology on weight truncation have a completely parallel version for local homology (\ie, derived completion)%
\footnote{For details about local homology, see \cite{Yeu20c}, whose notation is adopted here.}. Indeed, one can show that, in place of Lemma \ref{weight_subset_tor}, we have $\cD_{<w}(\GrA) \subset \cD_{\Imcomp}(\GrA)$ (see \cite[Proposition 2.42]{Yeu20c}). This allows one to develop formal analogues of Proposition \ref{RGam_Ce_weight_descend}, Theorem \ref{RGam_Ce_weight_ortho}, Theorem  \ref{three_term_SOD} and Theorem \ref{I_triv_weight_equiv}, where the functors $(\RGam_{I^+}, \Ce_{I^+}, \mathscr{L}_{[\geq w]}, \cL_{[\geq w]})$ are replaced by $(\LLam_{I^-}, \cE_{I^-}, \mathscr{R}_{\{\geq w\}}, \cR_{\{\geq w\}})$. 
However, this seems to be less useful for us because the results in the next section seems to have no analogue for this dual version.
\erm

For later use, we also show that the equivalence in Theorem \ref{I_triv_weight_equiv} have finite cohomological dimension. This is clear for the functor \eqref{I_triv_weight_functor_1} since it descends from an exact functor on the abelian categories. For the functor \eqref{I_triv_weight_functor_2}, we have the following
\blm  \label{I_triv_weight_cohom_dim}
The functor $\Ce_{I^+} \circ \mathscr{L}_{[\geq w]}  :  \cD(\Fgw) \ra \cD_{\IpTR}(\GrA)$ has finite cohomological dimension. 
\elm

\bpf
Since the exact functor $(-)^{\sharp}_{\geq w} : \GrA \ra \Mod(\Fgw)$ is essentially surjective, the statement in the Lemma is equivalent to the statement that the functor $\Ce_{I^+} \circ \cL_{[\geq w]} : \cD(\GrA) \ra \cD_{\IpTR}(\GrA)$ has finite cohomological dimension.
To show this, simply apply $\Ce_{I^+}$ to the exact triangle in the first row of \eqref{weight_trunc_exact_tri}. By Lemma \ref{weight_subset_tor}, we therefore have $ \Ce_{I^+}(\cL_{[\geq w]}(M)) \cong \Ce_{I^+}(M)$, which has cohomology in degrees $[p,q+r]$ if $M$ has cohomology in degrees $[p,q]$ (here $r$ is the number of generators of $I^+$). 
\epf

\section{Coherent subcategories}

The goal of this section is two-fold. We first show that the semi-orthogonal decomposition in Theorem \ref{three_term_SOD} restricts to a semi-orthogonal decomposition on the subcategory $\Dmcoh(\GrA)$ (see Theorem \ref{three_term_SOD_Dmcoh} below). Then we show that the equivalence in Theorem \ref{I_triv_weight_equiv} restricts to an equivalence between suitable bounded coherent subcategories (see Theorem \ref{I_triv_weight_equiv_Dsuitcoh}).

From now on, assume that the $\bZ$-graded ring $A$ is Noetherian. Then by Corollary \ref{Fgw_Noeth}, $\Fgw$ is also Noetherian, a fact that we will use without any more explicit mention.  We start with the following
\bpp  \label{weight_Cech_Dbcoh}
If $M \in \Dbcoh(\GrA)$, then we have $\Ce_{I^+}(M)^{\sharp}_{\geq w} \in \Dbcoh(\Fgw)$.
\epp

\bpf
It suffices to assume that $M$ is a finitely generated graded $A$-module concentrated in cohomological degree $0$. Recall from \eqref{Ce_i_RGam} that $H^p(\Ce_{I^+}(M))_i \cong H^p(X^+ , \widetilde{M(i)})$. Thus, by Lemma \ref{Fga_fg_module}, it suffices to show that 
$\bigoplus_{i \geq w} H^p(X^+ , \widetilde{M(i)})$ is finitely generated over $A_{\geq 0}$. 

Recall (see, e.g., \cite[Tag 0B5T]{Sta}) that if $\cF$ is a coherent sheaf on $X^+$ and $\cL$ is an ample invertible sheaf on $X^+$, then $\bigoplus_{i \geq w} H^p(X^+ , \cF \otimes \cL^{\otimes i})$ is finitely generated over $B_{\geq 0} := \bigoplus_{i \geq 0} H^0(X^+ , \cL^{\otimes i})$. 
In the present case, suppose that $A$ is positively $\tfrac{1}{d}$-Cartier, then our statement follows by applying this to $\cF = \widetilde{M}, \widetilde{M(1)},\ldots, \widetilde{M(d-1)}$ and $\cL = \widetilde{A(d)}$,
because the $\bN$-graded algebra $B_{\geq 0} := \bigoplus_{i \geq 0} H^0(X^+ , \widetilde{A(di)}$ is itself finite over $A^{(d)}$.
\epf

\brm  \label{weight_Cech_Dbcoh_remark}
Proposition \ref{weight_Cech_Dbcoh} (and its proof) is one of the major advantages of imposing the weight truncation. Namely, while $\Ce_{I^+}(M)$ has bounded cohomology, its cohomology groups $H^p(\Ce_{I^+}(M)) = \bigoplus_{i \in \bZ} H^p(X^+ , \tilde{M}(i))$ are \emph{not} finitely generated over $A$. 
This is because, while the sheaf $\bigoplus_{i \geq w} \widetilde{M(i)}$ is finitely generated over the sheaf $\pi^*(A_{\geq 0})$ of algebras, the sheaf $\bigoplus_{i \in \bZ} \widetilde{M(i)}$ is \emph{not} finitely generated over the sheaf $\pi^*(A)$ of algebras ({\it cf}. \cite[Tag 0897]{Sta}).
\erm

\bcor  \label{weight_Cech_Dsuitcoh}
For each $\spadesuit \in \{ \, \, , +,-,b\}$, 
if $M \in \Dsuit_{{\rm coh}}(\GrA)$, then we have $\Ce_{I^+}(M)^{\#}_{\geq w} \in \Dsuit_{{\rm coh}}(\Fgw)$.
\ecor

\bpf
Since $\Ce_{I^+}$ has bounded cohomological dimension, say $\Ce_{I^+}(\cD^{\leq p}(\GrA)) \subset \cD^{\leq p+m}(\GrA)$ and $\Ce_{I^+}(\cD^{\geq q}(\GrA)) \subset \cD^{\geq q}(\GrA)$, 
we have 
$H^p( \Ce_{I^+}(M)^{\sharp}_{\geq w} ) \cong H^p( \Ce_{I^+}(\tau_{\geq p-m}\tau_{\leq p}M)^{\sharp}_{\geq w} )$.
Apply Proposition \ref{weight_Cech_Dbcoh} to $\tau_{\geq p-m}\tau_{\leq p}M \in \Dbcoh(\GrA)$ to conclude that $H^p( \Ce_{I^+}(M)^{\sharp}_{\geq w} ) \in \Mod(\Fgw)$ is finitely generated.
\epf

\bdf
For each $\spadesuit \in \{ \, \, , +,-,b\}$, define the following subcategories of $\cD(\Fgw)$:
\begin{equation*}
\begin{split}
\Dsuit_{{\rm coh}, \, \IpTR}(\Fgw) \, &:= \, \Dsuit_{{\rm coh}}(\Fgw) \, \cap \, \cD_{\IpTR}(\Fgw) \\
\Dsuit_{{\rm coh}, \, \Torp}(\Fgw) \, &:= \, \Dsuit_{{\rm coh}}(\Fgw) \, \cap \, \cD_{\Torp}(\Fgw) 
\end{split}
\end{equation*}
\edf
Then Corollary \ref{weight_Cech_Dsuitcoh} gives the following
\bcor   \label{weight_Cech_Dsuitcoh_SOD}
For each $\spadesuit \in \{ \, \, , +,-,b\}$, the semi-orthogonal decomposition in Theorem \ref{RGam_Ce_weight_ortho}(4) restricts to a semi-orthogonal decomposition
\begin{equation}  \label{Dbcoh_Fgw_SOD}
\Dsuit_{{\rm coh}}(\Fgw) = \langle \, \Dsuit_{{\rm coh}, \, \IpTR}(\Fgw) \, , \, \Dsuit_{{\rm coh}, \, \Torp}(\Fgw) \, \rangle
\end{equation}
\ecor

We wish to combine this with the semi-orthogonal decomposition \eqref{weight_SOD} to obtain a three-term semi-orthogonal decomposition on the full subcategory $\Dmcoh(\GrA)$. To this end, we have to show that the weight truncation functor $\cL_{[\geq w]}$, and hence $\cL_{<w}$, preserve the full subcategory $\Dmcoh(\GrA) \subset \cD(\GrA)$. This follows from the following two lemmae: 

\blm  \label{preserve_Dmcoh_1}
If $\cM \in \Dmcoh(\Fgw)$, then $\mathscr{L}_{[\geq w]}(\cM) \in \Dmcoh(\GrA)$.
\elm

\bpf
Recall from Proposition \ref{Dmcoh_Dpc_right_Noeth} that $\Dbcoh(\Fgw) = \Dpc(\Fgw)$, so that $\cM$ may be represented by a bounded above complex $\cP^{\bullet}$ of free modules of finite rank. By \eqref{weight_tensor_free}, the complex $(\cP^{\bullet} \otimes_{\Fgw} \cF)^{\flat}$ is therefore in $\Dpc(\GrA) = \Dmcoh(\GrA)$.
\epf

\blm  \label{preserve_Dmcoh_2}
For each $\spadesuit \in \{ \, \, , +,-,b\}$, if $M \in \Dsuit_{{\rm coh}}(\GrA)$, then 
$M^{\sharp}_{\geq w} \in \Dsuit_{{\rm coh}}(\Fgw)$.
\elm

\bpf
By Lemma \ref{Fga_fg_module}, the exact functor $(-)^{\sharp}_{\geq w} : \Gr(A) \ra \Mod(\Fgw)$ sends finitely generated graded modules to finitely generated modules.
\epf

As a result, if we let 
\begin{equation*}
\begin{split}
\cD^-_{{\coh},  [\geq w]}(\GrA) \, &:= \, \Dmcoh(\GrA) \, \cap \, \cD_{[\geq w]}(\GrA) \\
\cD^-_{{\coh},  <w}(\GrA) \, &:= \, \Dmcoh(\GrA) \, \cap \, \cD_{<w}(\GrA) 
\end{split}
\end{equation*}
then we have the following
\bpp  \label{weight_SOD_Dmcoh}
The full subcategory $\cD^-_{{\coh},  [\geq w]}(\GrA) \subset \cD(\GrA)$  is the essential image of 
$\Dmcoh(\Fgw)$ under the fully faithful functor $\mathscr{L}_{[\geq w]} : \cD(\Fgw) \ra \cD(\GrA)$.
Moreover, the semi-orthogonal decomposition \eqref{weight_SOD} restricts to a semi-orthogonal decomposition
\begin{equation*}
\Dmcoh(\GrA) = \langle \, \cD^-_{\coh,< w}(\GrA) \, , \, \cD^-_{\coh,[\geq w]}(\GrA) \, \rangle
\end{equation*}
\epp

Combining Corollary \ref{weight_Cech_Dsuitcoh_SOD} and Proposition \ref{weight_SOD_Dmcoh}, we see that all the decomposition terms in \eqref{triple_SOD_terms} lie in $\Dmcoh(\GrA)$. As a result, we have the following
\bthm  \label{three_term_SOD_Dmcoh}
The semi-orthogonal decomposition in Theorem \ref{three_term_SOD} restricts to a semi-orthogonal decomposition 
\begin{equation*}
\Dmcoh(\GrA) \, = \, \langle \, \cD^-_{\coh,<w}(\GrA) \, , \, \mathscr{L}_{[\geq w]} ( \cD^-_{\coh, \, \IpTR}(\Fgw) )  \, , \,  \mathscr{L}_{[\geq w]} ( \cD^-_{\coh, \, \Torp}(\Fgw) ) \, \rangle
\end{equation*}
where the latter two semi-orthogonal components can be identified as
\begin{equation*}  
\begin{split}
\mathscr{L}_{[\geq w]} ( \cD^-_{\coh, \, \IpTR}(\Fgw) ) \, &= \, \{ \, M \in \cD^-_{\coh,[\geq w]}(\GrA) \, | \, \RGam_{I^+}(M) \in \cD_{<w}(\GrA) \, \} \\
\mathscr{L}_{[\geq w]} ( \cD^-_{\coh, \, \Torp}(\Fgw) )
\, &= \, \cD^-_{\coh, \Torp, [\geq w]}(\GrA) \, := \, \cD^-_{\coh, [\geq w]}(\GrA) \cap \cD_{\Torp}(\GrA)
\end{split}
\end{equation*}
\ethm

Next we show that the equivalence in Theorem \ref{I_triv_weight_equiv} restricts to an equivalence on coherent subcategories (see Theorem \ref{I_triv_weight_equiv_Dsuitcoh} below). 
Recall from Definition \ref{IpTR_coh_def} that, for each $\spadesuit \in \{ \, \, , +,-,b\}$, the full subcategory $\Dsuit_{{\rm coh}(\IpTR)}(\GrA) \subset \Dsuit_{\IpTR}(\GrA)$ is defined to be the essential image of $\Dsuit_{{\rm coh}}(\GrA)$ under the functor $\Ce_{I^+} : \Dsuit(\GrA) \ra \Dsuit_{\IpTR}(\GrA)$.
In view of Proposition \ref{coh_IpTR_QpGr}, this is the ``correct'' coherent subcategory to consider.

\bthm  \label{I_triv_weight_equiv_Dsuitcoh}
For each $\spadesuit \in \{ \, \, , +,-,b\}$, the equivalences in Theorem \ref{I_triv_weight_equiv} restricts to equivalences 
\begin{equation*}
\begin{tikzcd}
\Ce_{I^+} \circ \mathscr{L}_{[\geq w]} \, : \, \Dsuit_{{\rm coh}, \, \IpTR}(\Fgw) \ar[r, shift left = 2] \ar[r, phantom, "\simeq" description] 
& \Dsuit_{\coh(\IpTR)}(\GrA) \, : \, (-)^{\#}_{\geq w} \ar[l, shift left = 2]
\end{tikzcd}
\end{equation*}
\ethm

\bpf
The fact that the functor \eqref{I_triv_weight_functor_1} sends $\Dsuit_{\coh(\IpTR)}(\GrA)$ to 
$\Dsuit_{{\rm coh}, \, \IpTR}(\Fgw)$ is precisely the content of Corollary \ref{weight_Cech_Dsuitcoh}. 
For the other direction, notice that Lemma \ref{preserve_Dmcoh_1} establishes the statement for $\spadesuit = -$. The general case then follows from Lemma \ref{I_triv_weight_cohom_dim} by using a standard truncation argument as in the proof of Corollary \ref{weight_Cech_Dsuitcoh}.
\epf

\section{Boundedness}

In this section, we first give an alternative characterization of the full subcategory $\cD^-_{\coh,[\geq w]}(\GrA)$ 
(see, e.g., Proposition \ref{Dmcoh_geq_w_equiv} below).
This allows us to show that our three-term semi-orthogonal decomposition in Theorem \ref{three_term_SOD_Dmcoh} coincides with that of \cite{HL15} when the regularity assumptions in \cite{HL15} are satisfied (see Remark \ref{SOD_coincide_remark}).
We also follow the arguments of \cite{HL15} to give a sufficient condition for the semi-orthogonal decomposition in Theorem \ref{three_term_SOD_Dmcoh} to restrict to one on $\Dbcoh(\GrA)$ (see Theorem \ref{cond_then_reg} below). By using Lemma \ref{local_cohom_weight_bounded}, we are able to circumvent \cite[Proposition 3.31]{HL15} in the proof of an analogue of \cite[Lemma 3.36]{HL15}. This in turn allows us to weaken the ``Assumption (A)'' in \cite{HL15}.

For any Noetherian $\bZ$-graded ring $A$, let $\Gr_{<w}(A) \subset \Gr(A)$ be the Serre subcategory consisting of graded modules $M \in \Gr(A)$ such that $M_i = 0$ for all $i \geq w$. 
Let $\gr(A) \subset \Gr(A)$ be the full subcategory of finitely generated graded modules, and let $\gr_{<w}(A) := \gr(A) \cap \Gr_{<w}(A)$. 

Throughout this section, we fix a graded ideal $I^{\prime +} \subset A$ such that
\begin{equation}  \label{I_prime_plus_cond}
I^+ \subset I^{\prime +} \subset \sqrt{I^+}
\end{equation}
Since the notion of $I^{\infty}$-torsion modules, $I$-trivial complexes, etc, depend only on $\sqrt{I}$, all our previous discussions remain valid if we replace $I^+$ by $I^{\prime +}$ everywhere.

We start with the following simple

\blm
The Serre subcategory
$\gr_{<w}(A) \subset \gr(A)$ is the smallest Serre subcategory containing the essential image of $\gr_{<w}(A/I^{\prime +})$ 
under the (fully faithful) functor $\gr(A/I^{\prime +}) \ra \gr(A)$.
\elm  

\bpf
Any module in $\Gr_{<w}(A)$ is clearly $(I^+)^{\infty}$-torsion. Thus if $M \in \gr_{<w}$ then let $M[(I^{\prime +})^m] := \{ x \in M \, | \, (I^{\prime +})^m \cdot x = 0\}$, we have an increasing filtration 
$0 \subset M[I^{\prime +}] \subset M[(I^{\prime +})^2] \subset \ldots$ whose union is $M$. Since $M$ is Noetherian, this must stabilize after finitely many terms. Since all the successive quotients in this filtration lies in $\gr_{<w}(A/I^{\prime +})$, we have the desired result.
\epf

\bcor  \label{Dbcoh_a_ess_im}
The triangulated subcategory $\cD^b_{\coh,<w}(\GrA) \subset \Dbcoh(\GrA)$ is the smallest triangulated subcategory of $\Dbcoh(\GrA)$ that contains the essential image of $\cD^b_{\coh,<w}(\Gr(A/I^{\prime +}))$ under the functor $\Dbcoh(\Gr(A/I^{\prime +})) \ra \Dbcoh(\GrA)$.
\ecor

A completely parallel proof also shows the following
\blm  \label{Torp_ess_im}
The Serre subcategory $\gr(A) \cap \Torp(A) \subset \gr(A)$ is the smallest Serre subcategory containing the essential image of $\gr(A/I^{\prime +}) \ra \gr(A)$.
Therefore, the triangulated subcategory $\cD^b_{\coh,\Torp}(\GrA) \subset \Dbcoh(\GrA)$ is the smallest triangulated subcategory that contains the essential image of the functor $\Dbcoh(\Gr(A/I^{\prime +})) \ra \Dbcoh(\GrA)$.
\elm

\bpp  \label{D_geq_a_coh_prop1}
For any $M \in \cD^-_{\coh}(\GrA)$, we have 
$M \in \cD^-_{\coh,[\geq w]}(\GrA)$ if and only if $M \otimes_A^{{\bm L}} (A/I^{\prime +}) \in \cD^-_{\coh,[\geq w]}(\Gr (A/I^{\prime +}))$.
\epp

\bpf
The implication ``$\Rightarrow$'' is clear. For the converse, suppose that $M \otimes_A^{{\bm L}} (A/I^{\prime +}) \in \cD^-_{[\geq w],\coh}(\Gr (A/I^{\prime +}))$. Since the functor $-\otimes_A^{{\bm L}} (A/I^{\prime +}) $ is left adjoint to the restriction of scalar functor $(-)_A$, we have $\RHom_A(M,Q_A) \simeq 0$ for each $Q \in \cD^-_{\coh,<w}(\Gr(A/I^{\prime +}))$. 
In view of Corollary \ref{Dbcoh_a_ess_im}, we have in particular $\RHom_A(M,K) \simeq 0$ for each $K \in \cD^b_{\coh,<w}(\GrA)$.
If $M$ is not in $\cD^-_{[\geq w],\coh}(\GrA)$, then by Proposition \ref{weight_SOD_Dmcoh}, there exists a nonzero map to some $N \in \cD^-_{\coh,<w}(\GrA)$, which must remain nonzero after passing to the truncation $N \ra \tau_{\geq m}(N)$ for some $m \in \bZ$. Since $\tau_{\geq m}(N) \in \cD^b_{\coh,<w}(\GrA)$, this gives a contradiction.
\epf

Now suppose that $B$ is a Noetherian ($-\bN$)-graded ring. \ie, $B$ is a $\bZ$-graded ring such that $B_i = 0$ for all $i >0$. Then we have $I^-(B) = B_{<0}$, and hence $B/I^-(B) = B_0$. In this case, weight truncation can often be performed inductively:
\blm  \label{B_weight_trunc}
Suppose that $M \in \cD_{<w}(\Gr(B))$, then we have 
$\cL_{[\geq w-1]}(M) \cong M_{w-1}(-w+1) \otimes_{B_0}^{{\bm L}} B$.
\elm

\bpf 
Take the adjunction counit $M_{w-1}(-w+1) \otimes_{B_0}^{{\bm L}} B \ra M$, which is clearly a quasi-isomorphism in weight $a-1$, so that its cone lies in $\cD_{<w-1}(\Gr(B))$, and hence is equal to $\cL_{<w-1}(M)$, since we have $M_{w-1}(-w+1) \otimes_{B_0}^{{\bm L}} B \in \cD_{[\geq w-1]}(\Gr(B))$.
\epf

The functor $-\otimes^{\bm L}_{B} B_0$ preserves both the subcategories $\cD_{[\geq w]}$ and $\cD_{<w}$:

\blm  \label{D_geq_a_lem1}
If $M \in \cD_{[\geq w]}(\Gr(B))$ then  $M \otimes^{{\bm L}}_{B} B_0 \in \cD_{[\geq w]}(\Gr(B_0))$.

If $M \in \cD_{<w}(\Gr(B))$ then  $M \otimes^{{\bm L}}_{B} B_0 \in \cD_{<w}(\Gr(B_0))$.
\elm

\bpf
The first statement is obvious. For the second statement, we apply Lemma \ref{B_weight_trunc}.
It is clear that $(M_{w-1}(-w+1) \otimes_{B_0}^{{\bm L}} B) \otimes^{{\bm L}}_{B} B_0$ is concentrated in weight $w-1$.
Thus, a repeated application of Lemma \ref{B_weight_trunc} gives a sequence of maps
\begin{equation*}
M = \cL_{<w}(M) \raq \cL_{<w-1}(M) \raq \cL_{<w-2}(M) \raq \ldots
\end{equation*}
such that $\cone[\, \cL_{<i+1}(M) \ra \cL_{<i}(M) \,] \otimes^{{\bm L}}_{B} B_0$ is concentrated in weight $i$. Thus, for any $i < w$, the weight truncation  
$\cL_{[\geq i]}(M) = \cone(M \ra \cL_{<i}(M))[-1]$ satisfies 
\begin{equation*}
\parbox{40em}{$\cL_{[\geq i]}(M) \otimes_B^{{\bm L}} B_0$ is concentrated in weight $[i,w-1]$. }
\end{equation*} 
Since the sequence of maps
\begin{equation*}
\cL_{[\geq w-1]}(M) \ra \cL_{[\geq w-2]}(M) \ra \ldots \ra M
\end{equation*}
exhibits $M$ as a homotopy colimit in $\cD(\Gr(B))$, and since homotopy colimit commutes with the functor $-\otimes_B^{{\bm L}} B_0$, we have 
$M \otimes^{{\bm L}}_{B} B_0 \in \cD_{<w}(\Gr(B_0))$.
\epf

We also have the following

\blm  \label{D_geq_a_lem2}
If $M \in \Dmcoh(\Gr(B))$ is a nonzero object, then $M \otimes^{{\bm L}}_{B} B_0 \in \Dmcoh(\Gr(B_0))$ is also nonzero.
\elm

\bpf
Take the highest nonvanishing cohomology degree $H^p(M) \neq 0$. Then by the Nakayama lemma for $\bN$-graded rings, we have $0 \neq H^p(M) \otimes_B B_0 = H^p( M \otimes^{{\bm L}}_{B} B_0)$.
\epf

Combining these two lemmae, we have

\bpp  \label{B_geq_w_B_zero}
For any $N \in \Dmcoh(\Gr(B))$, we have
$N \in \cD^-_{\coh,[\geq w]}(\Gr(B))$ if and only if $N \otimes^{{\bm L}}_B B_0 \in \cD^-_{\coh,[\geq w]}(\Gr(B_0))$.
\epp

\bpf
The direction ``$\Rightarrow$'' is obvious. For the direction ``$\Leftarrow$'', suppose that $N \notin \cD^-_{\coh,[\geq w]}(\Gr(B))$ so that $\cL_{<w}(N) \neq 0$. Then by Lemma \ref{D_geq_a_lem1}, the exact triangle 
\begin{equation*}
\ldots \raq \cL_{[\geq w]}(N) \otimes^{{\bm L}}_B B_0 \raq N \otimes^{{\bm L}}_B B_0 \raq \cL_{[< w]}(N) \otimes^{{\bm L}}_B B_0 \xraq{[1]} \ldots
\end{equation*}
is precisely the weight truncation sequence for $N \otimes^{{\bm L}}_B B_0$ in $\cD(\Gr(B_0))$. By Lemma \ref{D_geq_a_lem2}, we have $ \cL_{<w}(N) \otimes^{{\bm L}}_B B_0 \neq 0$, which therefore shows that $N \otimes^{{\bm L}}_B B_0 \notin \cD^-_{\coh,[\geq w]}(\Gr(B_0))$.
\epf

\bpp  \label{Dmcoh_geq_w_equiv}
For any $M \in \Dmcoh(\GrA)$, 
the followings are equivalent:
\begin{enumerate}
	\item $M \in \cD^-_{\coh,[\geq w]}(\GrA)$
	\item $M \otimes^{{\bm L}}_A (A/(I^- + I^+)) \in \cD^-_{\coh,[\geq w]}(\Gr(A/(I^- + I^+))$. 
	\item $M \otimes^{{\bm L}}_A (A/\sqrt{I^- + I^+}) \in \cD^-_{\coh,[\geq w]}(\Gr(A/\sqrt{I^- + I^+})$. 
\end{enumerate}
\epp

\bpf
Take $B = A/I^{\prime +}$ for $I^{\prime +} = I^+$ or $I^{\prime +} = \sqrt{I^+}$.
In the former case, we have $B_0 = B/I^-(B) = A/(I^- + I^+)$. 
In the latter case, $B_0$ is a subring of the reduced ring $B$, and hence is reduced. 
In other words, $I^- + \sqrt{I^+} \subset A $ is equal to its radical, and must therefore be equal to $\sqrt{I^- + I^+}$. Thus, we have $B_0 = A/\sqrt{I^- + I^+}$, and it suffices to show that 
\begin{equation*}
M \in \cD^-_{\coh,[\geq w]}(\GrA) \quad \Leftrightarrow \quad 
M \otimes^{{\bm L}}_A B_0 \in \cD^-_{\coh,[\geq w]}(\Gr(B_0))
\end{equation*}
for $B = A/I^{\prime +}$, where $I^{\prime +}$ is any graded ideal satisfying \eqref{I_prime_plus_cond}.
Take $N := M \otimes_A^{{\bm L}} B \in  \Dmcoh(\Gr(B))$.
The result then follows from Propositions \ref{D_geq_a_coh_prop1} and \ref{B_geq_w_B_zero}.
\epf

Now we give a sufficient condition for weight truncation to preserve $\Dbcoh(\GrA)$ (see Theorem \ref{cond_then_reg}). The arguments for Lemma \ref{B_flat_then_reg}, Propositon \ref{Torp_reg} and Theorem \ref{cond_then_reg} below are adapted from those in \cite{HL15}. However, we weaken the assumption $(A)$ in {\it loc.cit.}.
	
Take a graded ideal $I^{\prime +} \subset A$ satisfying \eqref{I_prime_plus_cond}. 
Consider the conditions
\begin{equation} \label{two_conds}
\parbox{40em}{(a) The graded ring $B := A/I^{\prime +}$ has finite Tor-dimension over the subring $B_0 \subset B$. \\
(b) As a quotient, the ring $B_0 = B/I^-(B)$ has finite Tor-dimension over $B$.\\
(c) $A/I'^+ \in \cD_{[\geq 0]}(\GrA)$.}
%
\end{equation}
Notice that by cocontinuity, condition (c) is equivalent to the condition that the restriction of scalar functor sends $\cD_{[\geq w]}(\Gr(A/I'^+))$ to $\cD_{[\geq w]}(\GrA)$.

Under the first two conditions, we have the following
\blm  \label{B_flat_then_reg}
Suppose that \eqref{two_conds}(a) holds, then weight truncation for $B$ preserves $\Dbcoh(\Gr(B))$. \ie, for all $M \in \Dbcoh(\Gr(B))$, we have $\cL_{< w}(M) \in \Dbcoh(\Gr(B))$. 

Suppose that \eqref{two_conds}(b) holds, then for every $M \in \Dbcoh(\Gr(B))$ there exists $i \in \bZ$ such that $M \in \cD^b_{\coh,[\geq i]}(\Gr(B))$.
\elm

\bpf
Since $B$ is ($-\bN$)-graded, there exists some $w' \in \bZ$ such that $M \in \cD_{<w'}(\Gr(B))$. Apply Lemma \ref{B_weight_trunc}, we see that $\cL_{[\geq w'-1]}(M) \cong M_{w-1}(-w'+1) \otimes_{B_0}^{{\bm L}} B$, which is in $\Dbcoh(\Gr(B))$ by the assumption \eqref{two_conds}(a). Thus, $\cL_{< w'-1}(M) \in \Dbcoh(\GrA)$. A repeated application of the argument then shows that $\cL_{<w}(M) \in \Dbcoh(\GrA)$ for all $w \in \bZ$.

For the second statement, the assumption \eqref{two_conds}(b) guarantees that $M \otimes_B^{{\bm L}} B_0 \in \Dbcoh(\Gr(B_0))$. Since $B_0$ is concentrated in weight $0$, any finitely generated graded $B_0$-module must be concentrated in finitely many weight components. Hence, $M \otimes_B^{{\bm L}} B_0 \in \cD_{[\geq i]}(\Gr(B_0))$ for some $i \in \bZ$.
By Proposition \ref{B_geq_w_B_zero}, this is precisely the sought for statement.
\epf

\bpp  \label{Torp_reg}
Suppose that conditions \eqref{two_conds}(a) and \eqref{two_conds}(c) hold, then 
\begin{enumerate}
	\item $\cD^b_{\coh,\Torp,[\geq w]}(\GrA)$ is the smallest triangulated subcategory containing the essential image of the functor
	$\cD^b_{\coh,[\geq w]}(\Gr(A/I^{\prime +})) \ra \Dbcoh(\GrA)$.
	\item For any $M \in \cD^b_{\coh, \Torp}(\GrA)$, we have $\cL_{[\geq w]}(M) \in \Dbcoh(M)$.
\end{enumerate} 
If \eqref{two_conds}(b) also hold, then we also have
\begin{enumerate}[resume]
	\item for any $M \in \cD^b_{\coh,\Torp}(\GrA)$ there exists $i \in \bZ$ such that $M \in \cD^b_{\coh,[\geq i]}(\Gr(A))$.
\end{enumerate}
\epp

\bpf
Consider the following full subcategories of $\Dbcoh(\Gr(B))$:
\begin{equation*}
\begin{split}
\cE_1 \,& := \, {\rm EssIm}( \, \cD^b_{\coh,<w}(\Gr(B)) \raq \Dbcoh(\GrA) \, ) \\
\cE_2 \,& := \, {\rm EssIm}( \, \cD^b_{\coh,[\geq w]}(\Gr(B)) \raq \Dbcoh(\GrA) \, )
\end{split}
\end{equation*} 

We shall adopt the notation from \cite[Appendix A]{Yeu20c}.
In particular, for any full subcategory $\cE \subset \Dbcoh(\GrA)$, we denote by $\langle \cE \rangle$ the smallest triangulated subcategory containing $\cE$. In this notation, Corollary \ref{Dbcoh_a_ess_im} asserts that $\langle \cE_1 \rangle = \cD^b_{\coh,<w}(\GrA)$. 
Condition \eqref{two_conds}(c) implies that $\cE_2 \subset \cD^b_{\coh, \Torp, [\geq w]}(\GrA)$, so that $\cE_1$ and $\cE_2$ are strongly orthogonal, \ie, $\Hom_{\cD(\GrA)}(E_2,E_1[i]) = 0$ for all $i \in \bZ$, $E_1 \in \cE_1$ and $E_2 \in \cE_2$.
Recall from \cite[Corollary A.9]{Yeu20c} that this implies that $\langle \cE_1 \rangle \vec{*} \langle \cE_2 \rangle = \langle \cE_1  \vec{*} \, \cE_2 \rangle$.
By Lemma \ref{B_flat_then_reg},  we have $\Dbcoh(\Gr(B)) = \langle \cD^b_{\coh,<w}(\Gr(B)), \cD^b_{\coh,[\geq w]}(\Gr(B)) \rangle$ under condition \eqref{two_conds}(a), so that 
\begin{equation*}
\cE \, := \, {\rm EssIm}( \, \cD^b_{\coh}(\Gr(B)) \raq \Dbcoh(\GrA) \, ) \, \subset \,  \cE_1 * \cE_2 
\end{equation*}
By Lemma \ref{Torp_ess_im}, we have $\langle \cE \rangle = \cD^b_{\coh,\Torp}(\GrA)$. Combining these facts, we have
\begin{equation*}
\begin{split}
\langle \cE \rangle \, &\subset \, \langle \cE_1 \vec{*} \, \cE_2 \rangle  \, = \, \langle \cE_1 \rangle \vec{*} \langle \cE_2 \rangle \\
 &\subset \,  (\cD^b_{\coh,<w}(\GrA)) \, \vec{*} \, (\cD^b_{\coh, \Torp, [\geq w]}(\GrA))  \\
 &\subset \,  \cD^b_{\coh,\Torp}(\GrA) \, = \, \langle \cE \rangle
\end{split}
\end{equation*}
Since the first and last term are the same, we must have equalities. 
In particular, the equality for the second inclusion implies that $\langle \cE_2 \rangle =  \cD^b_{\coh, \Torp, [\geq w]}(\GrA)$, which is the first sought for statement. The equality for the third inclusion is precisely the second sought for statement.

If \eqref{two_conds}(b) holds, then applying the second statement of Lemma \ref{B_flat_then_reg}, together with \eqref{two_conds}(c), we see that for every object $N \in \cE$ there is some $i \in \bZ$ such that $N \in \cD^b_{\coh,\Torp,[\geq i]}(\GrA)$. Since $\cD^b_{\coh,\Torp}(\GrA) \, = \, \langle \cE \rangle$, every $M \in \cD^b_{\coh,\Torp}(\GrA)$ also has this property.
\epf

\bthm  \label{cond_then_reg}
Suppose that the conditions \eqref{two_conds}(a),(b),(c) hold. Then
\begin{enumerate}
	\item For every object $M \in \Dbcoh(\GrA)$, there exists some $i \in \bZ$ such that $M \in \cD^b_{\coh,[\geq i]}(\GrA)$.  
	\item Weight truncation for  $A$ preserves $\Dbcoh(\GrA)$. \ie, for any $M \in \Dbcoh(\GrA)$, we have $\cL_{<w}(M) \in \Dbcoh(\GrA)$.
\end{enumerate}
\ethm

\bpf
Let $f_1,\ldots,f_r \in A$ be a set of elements of positive degrees $\deg(f_i) = d_i > 0$ that generate $I^+$, and let $
K^{\bullet}(A,f_1,\ldots,f_r) \, = \, \bigwedge_A (A \theta_1 \oplus \ldots \oplus A \theta_r)
$ 
be the (cohomological) Koszul complex, which is a finite complex of free graded $A$-modules with a set $\{\wedge_{s \in S} \,  \theta_s\}_{S \subset \{1,\ldots,r\}}$ of $2^r$ generators of weight $-\sum_{s \in S} d_s$ and cohomological degree $|S|$. Moreover, the differentials of the Koszul complex satisfies 
\begin{equation}  \label{Koszul_diff_in_I}
d(K^{\bullet}(A,f_1,\ldots,f_r)) \, \subset \, I^+ \cdot K^{\bullet}(A,f_1,\ldots,f_r)
\end{equation} 
For any $M \in \Dbcoh(\GrA)$, let $K^{\bullet}(M,f_1,\ldots,f_r) := K^{\bullet}(A,f_1,\ldots,f_r) \otimes_A M$. Then \eqref{Koszul_diff_in_I} implies that \begin{equation*}
K^{\bullet}(M,f_1,\ldots,f_r) \otimes_{A}^{{\bm L}} B \, \, \, \cong \, \bigoplus_{S \subset \{1,\ldots,r\}} (M \otimes_A^{{\bm L}} B) \,( \, \textstyle \sum_{s \in S} d_s \,)[ \,- |S| \, ]
\end{equation*}
where $B := A/I^{\prime +}$.
By Proposition \ref{D_geq_a_coh_prop1}, we therefore see that $K^{\bullet}(M,f_1,\ldots,f_r) \in \cD^b_{\coh,[\geq i]}(\GrA)$ if and only if $M \in \cD^b_{\coh,[\geq i]}(\GrA)$.
Since we always have $K^{\bullet}(M,f_1,\ldots,f_r) \in \cD^b_{\coh,\Torp}(\GrA)$, the first statement of the present Theorem follows from Proposition \ref{Torp_reg}(3).

Let $K_{\bullet}(A,f_1^j,\ldots,f_r^j)$ be the homological Koszul complex, \ie, it is the $A$-linear dual $\Homcom_A(-,A)$ of $K^{\bullet}(A,f_1^j,\ldots,f_r^j)$. Thus, it is a finite complex of free graded $A$-modules with a set $\{\wedge_{s \in S} \,  \theta_s^{\vee}\}_{S \subset \{1,\ldots,r\}}$ of $2^r$ generators of weight $(\sum_{s \in S} d_s)j$ and cohomological degree $-|S|$.
Let $K_{\bullet}(M,f_1^j,\ldots,f_r^j) := M \otimes_A  K_{\bullet}(A,f_1^j,\ldots,f_r^j)$, then by statement (1) we have just proved, we have 
\begin{equation*}
\cone \, [ \, M \raq K_{\bullet}(M,f_1^j,\ldots,f_r^j) \, ] \, \in \, \cD_{[\geq w]}(\GrA)
\qquad \text{for } \, j \gg 0 
\end{equation*}
As a result, we have 
$\cL_{<w}(M) \cong \cL_{<w} ( K_{\bullet}(M,f_1^j,\ldots,f_r^j))$ for $j \gg 0$.
Since $K_{\bullet}(M,f_1^j,\ldots,f_r^j) \in \cD^b_{\coh,\Torp}(\GrA)$, the second statement follows from Proposition \ref{Torp_reg}(2).
\epf

Notice that if weight truncation preserves $\Dbcoh(\GrA)$,
then the semi-orthogonal decomposition in Theorem \ref{three_term_SOD_Dmcoh} restricts to a semi-orthogonal decomposition 
\begin{equation} \label{three_term_SOD_Dbcoh}
\Dbcoh(\GrA) \, = \, \langle \, \cD^b_{\coh,<w}(\GrA) \, , \, \mathscr{L}_{[\geq w]} ( \cD^b_{\coh, \, \IpTR}(\Fgw) )  \, , \,  \cD^b_{\coh, \Torp, [\geq w]}(\GrA) \, \rangle
\end{equation}
Moreover, the latter two semi-orthogonal components can be described as in Theorem \ref{three_term_SOD_Dmcoh}, with $\cD^-$ replaced by $\cD^b$ everywhere.


The following is the main class of examples of Noetherian $\bZ$-graded rings that satisfies the conditions \eqref{two_conds}(a),(b),(c):

\bpp  \label{smooth_then_cond}
If $A$ is a  $\bZ$-graded ring finitely generated over a field $k$ of characteristic zero, and if the underlying ungraded algebra $A$ is smooth over $k$, then we have
\begin{enumerate}
	\item The algebras $B^+ = A/I^+$, $B^- = A / I^-$ and $B_0 = A/(I^- + I^+)$ are smooth over $k$. 
	\item The projections $\rho^+ : \Spec \, B^+ \ra \Spec \, B_0$ and $\rho^- : \Spec \, B^- \ra \Spec \, B_0$ are locally trivial bundle of weighted affine spaces.
	\item Along each connected component $Z_i \subset \Spec\, B_0$, we have $\dim((\rho^+)^{-1}(Z_i)) + \dim((\rho^-)^{-1}(Z_i)) = \dim(Z_i) + \dim(A)$.
\end{enumerate}
Therefore, the conditions  \eqref{two_conds}(a),(b),(c) are satisfied.
\epp

\bpf
Notice that $B_0 = A/(I^- + I^+)$ is the weight zero part of both $B^+ = A/I^+$ and $B^- = A/I^-$, so that the second statement make sense. 
These statements are then a special case of a result of Bia\l{}ynicki-Birula \cite{BB73}. 
The conditions \eqref{two_conds}(a),(b),(c) then follows (see, \eg, \cite[Lemma 2.7]{HL15}).
\epf

As a consequence, we have the following
\bthm  \label{smooth_then_reg_and_perf}
If $A$ is a  $\bZ$-graded ring finitely generated over a field $k$ of characteristic zero, and if the underlying ungraded algebra $A$ is smooth over $k$, then the weight truncation functor $\cL_{[\geq w]}$ preserves $\Dbcoh(\GrA) = \Dperf(\GrA)$.
\ethm


For applications in Section \ref{Dcat_flip_flop_sec} below, we will be mostly interested in the case when the semi-orthogonal decomposition \eqref{Dbcoh_Fgw_SOD} restricts to one on $\Dperf(\Fgw)$. 
\blm
If the weight truncation functor $\cL_{[\geq w]}$ preserves $\Dperf(\GrA)$, then the semi-orthogonal decomposition \eqref{Dbcoh_Fgw_SOD} restricts to one on $\Dperf(\Fgw)$.
\elm

\bpf
Recall that there exists some $n_0 \geq 0$ such that any sequence $\Ce_{I^+}(A)(i), \ldots, \Ce_{I^+}(A)(i+n_0)$ of length $n_0 + 1$ is a set of compact generators for $\cD_{\IpTR}(\GrA)$. 
Take some $c \in \bZ$ such that $\RGam_{I^+}(A) \in \cD_{<c}(\GrA)$ (see Lemma \ref{local_cohom_weight_bounded}), so that $\Ce_{I^+}(A)(i)^{\sharp}_{\geq w} = A(i)^{\sharp}_{\geq w}$ for any $i \geq c - w$. Thus, for these values of $i$, we have $\cL_{[\geq w]}(\Ce_{I^+}(A)(i)) = \cL_{[\geq w]}( A(i) )$, which is in $\Dperf(\GrA)$ by assumption.
Since $\cD_{\IpTR}(\GrA)_c$ is split generated by these objects, we see that $\cL_{[\geq w]}(C) \in \Dperf(\GrA)$ for all $C \in \cD_{\IpTR}(\GrA)_c$.
In view of Lemma \ref{D_geq_w_gen}, this shows that the functor $\Ce_{I^+,\geq w}$ in \eqref{RGam_Ce_weight} preserves $\Dperf(\Fgw)$.
\epf

\brm  \label{SOD_coincide_remark}
Proposition \ref{Dmcoh_geq_w_equiv} allows us to compare our construction of weight truncation with the ones in \cite{HL15} and \cite{BFK19}. Namely, if the assumptions (L+) and (A) in \cite{HL15} are satisfied, then so does our assumptions \eqref{two_conds}. 
By comparing  Proposition \ref{Dmcoh_geq_w_equiv}  with \cite[Definition 2.8]{HL15} for $\mathfrak{X} := [\Spec \, A / \bG_m]$, we see that 
\begin{equation*}
\cD^b_{\coh,\Torp,[\geq w]}(\GrA) ) \, = \, \cD^b_{\mathfrak{X}^u}(\mathfrak{X})_{\geq w}
\end{equation*}
Then, notice that Lemma \ref{Torp_SOD} restricts to $\cD^b_{\coh,\Torp}(\GrA)$. Comparing this with \cite[Theorem 2.10(5)]{HL15}, we see that 
\begin{equation*}
\cD^b_{\coh,<w}(\GrA) \, = \, \cD^b_{\mathfrak{X}^u}(\mathfrak{X})_{< w}
\end{equation*}
Finally, if we compare the three term semi-orthogonal decomposition in \eqref{three_term_SOD_Dbcoh} with the corresponding one in \cite[Theorem 2.10(6)]{HL15}, then we see that 
\begin{equation*}
\mathscr{L}_{[\geq w]} ( \cD^b_{\coh, \, \IpTR}(\Fgw) ) \, = \, \mathbf{G}_w
\end{equation*}
This shows that the three-term semi-orthogonal decomposition of \cite{HL15}, and hence of \cite{BFK19}, coincides with the one in \ref{three_term_SOD_Dbcoh} in this abelian case. 
\erm





\section{The case of non-affine base}
In this section, we provide the formal arguments to extend our previous discussion to the case of non-affine base. More precisely, we work in the following setting:
\begin{equation}  \label{sheaf_A_setting}
\parbox{40em}{$Y$ is a Noetherian separated scheme, and $\cA$ is a quasi-coherent sheaf of Noetherian $\bZ$-graded rings on $Y$, such that $\cA_0$ (and hence every $\cA_i$) is coherent over $\cO_X$.}
\end{equation}

Denote by $\Gr(\cA)$ the category of quasi-coherent graded $\cA$-modules. Then for any quasi-coherent sheaf of graded ideal $\scI \subset \cA$, there is a semi-orthogonal decomposition 
\begin{equation*}  
\cD(\Gr(\cA)) \, = \, \langle \, \cD_{\cITR}(\Gr(\cA)) \, , \,  \cD_{\cIIT}(\Gr(\cA)) \, \rangle
\end{equation*}
whose restriction to each open affine subscheme is precisely \eqref{local_cohom_SOD}. More precisely, if we take the decomposition triangle 
\begin{equation}  \label{RGam_Ce_seq_sheaf}
\ldots \raq \RGam_{\scI}(\cM) \xraq{\epsilon_M} \cM \xraq{\eta_{\cM}} \Ce_{\scI}(\cM) \xraq{\delta_{\cM}} \RGam_{\scI}(\cM)[1] \raq \ldots 
\end{equation}
for $\cM\in \cD(\Gr(\cA))$, then the value on  $U = \Spec \, R \subset Y$ is precisely the \eqref{RGam_Ce_seq} for $M = \cM(U)$ (see, \eg, \cite[Section 5]{Yeu20c} for details). Once we have defined these functors, their properties can then be checked affine-locally.

One can also consider weight truncation for pairs $(Y,\cA)$ in \eqref{sheaf_A_setting}. 
As in the above discussion of local cohomology, it suffices to construct the relevant weight truncation functors that reduce to the ones above over any open affine subscheme $\Spec \, R \subset Y$. Then the properties of such functors can be checked locally.
We start with the following
\bdf
Let $\cD_{[\geq w]}(\Gr(\cA)) \subset \cD(\Gr(\cA))$ be the smallest strictly full triangulated subcategory  closed under small coproducts, and containing the object of the form
\begin{equation}  \label{gen_obj_Dgw}
\cF \otimes_{\cO_Y}^{{\bm L}} \cA(-i) \, , \qquad \text{where} \quad  \cF \in \cD_{\perf}(\QCoh(Y)) \quad \text{and} \quad  i \geq w
\end{equation}
\edf

Clearly, each of the objects of the form \eqref{gen_obj_Dgw} is compact in $\cD(\Gr(\cA))$, and hence also in $\cD_{[\geq w]}(\Gr(\cA))$. Thus, $\cD_{[\geq w]}(\Gr(\cA))$ is compactly generated%
, and the inclusion functor $\cD_{[\geq w]}(\Gr(\cA)) \rinto \cD(\Gr(\cA))$ preserves small coproducts.
By the Brown-Neeman representability theorem \cite[Theorem 4.1]{Nee96}, this inclusion therefore has a right adjoint, which will be denoted as
\begin{equation}
\cL_{[\geq w]} \, : \, \cD(\Gr(\cA)) \raq \cD_{[\geq w]}(\Gr(\cA))
\end{equation}

Since the full triangulated subcategory $\cD_{[\geq w]}(\Gr(\cA)) \subset \cD(\Gr(\cA))$ is right admissible, there exists a semi-orthogonal decomposition of the form
\begin{equation}  \label{weight_SOD_nonaffine}
\cD(\Gr(\cA)) \, = \, \langle \, \cD_{<w}(\Gr(\cA)) \, , \, \cD_{[\geq w]}(\Gr(\cA))  \, \rangle
\end{equation}
where $\cD_{<w}(\Gr(\cA)) := \cD_{[\geq w]}(\Gr(\cA))^{\perp}$. 
Alternatively, it can be characterized as follows:
\blm  \label{Dlw_nonaffine}
An object $\cM \in \cD(\Gr(\cA))$ is in $\cD_{<w}(\Gr(\cA))$ if and only if its $i$-th weight component $\cM_i \in \cD(\QCoh(Y))$ is zero for all $i \geq w$.
\elm

\bpf
For any $\cM \in \cD(\Gr(\cA))$, we have $\cM \in \cD_{[\geq w]}(\Gr(\cA))^{\perp} $ if and only if 
\begin{equation}  \label{left_orth_gen}
\Hom_{\cD(\Gr(\cA))}( \cF \otimes_{\cO_Y}^{{\bm L}} \cA(-i) , \cM[j]) = 0 \quad \text{for all} \quad \cF \in \cD_{\perf}(\Gr(\cA)) , \,  i \geq w \, \, \text{and} \, \, j \in \bZ 
\end{equation}
Indeed, by the simple fact \eqref{left_orth_coprod} below, applied to $\cD := \cD(\Gr(\cA))$ and $X := \cM$, we see that the objects \eqref{gen_obj_Dgw} are in $\,^{\perp}({\bm \Sigma} \cM)$ if and only if $\cD_{[\geq w]}(\Gr(\cA)) \subset \,^{\perp}({\bm \Sigma} \cM)$.
\begin{equation} \label{left_orth_coprod}
\parbox{40em}{Suppose $\cD$ is a triangulated category that admits small coproducts. Then for any $X \in \cD$, the full subcategory $\! \,^{\perp}({\bm \Sigma} X) := \{ \, Y \in \cD \, | \, \Hom_{\cD}(Y,X[i]) = 0 \text{ for all } i \in \bZ \, \}$ is a triangulated subcategory that is closed under small coproducts.}
\end{equation}

Notice that we have
\begin{equation*}
\Hom_{\cD(\Gr(\cA))}( \cF \otimes_{\cO_Y}^{{\bm L}} \cA(-i) , \cM[j]) \, \cong \, \Hom_{\cD(\QCoh(Y))} ( \cF , \cM_i )
\end{equation*}
Recall that $\cD(\QCoh(Y))$ is compactly generated by $\cD_{\perf}(\QCoh(Y))$ (see, e.g., \cite[Corollary 2.3, Proposition 2.5]{Nee96} and \cite[Theorem 3.1.1]{BVdB03}). The result therefore follows from the characterization \eqref{left_orth_gen} of $\cD_{[\geq w]}(\Gr(\cA))^{\perp}$.
\epf

By the semi-orthogonal decomposition \eqref{weight_SOD_nonaffine}, we see that the inclusion $\cD_{<w}(\Gr(\cA)) \rinto \cD(\Gr(\cA))$ has a left adjoint, which we denote as
\begin{equation}
 \cL_{< w} \, : \, \cD(\Gr(\cA)) \raq \cD_{< w}(\Gr(\cA))
\end{equation}

For any $\cM \in \cD(\Gr(\cA))$, the semi-orthogonal decomposition \eqref{weight_SOD_nonaffine} then gives us a decomposition sequence
\begin{equation}  \label{weight_decomp_seq_nonaffine}
\ldots \raq \cL_{[\geq w]}(\cM) \raq \cM \raq \cL_{< w}(\cM) \raq  \cL_{[\geq w]}(\cM)[1] \raq \ldots
\end{equation}

Given any open affine subscheme $U = \Spec \, R \subset Y$, let $A := \cA(U)$. If $\cM$ is of the form \eqref{gen_obj_Dgw}, then its restriction to $U$ has the form 
$K \otimes_R^{{\bm L}} A (-i)$, where $\cF \in \cD_{\perf}(R)$ and $i \geq w$.
Since $R$ split generates  $\cD_{\perf}(R)$, we see that these are all contained in $\cD_{[\geq w]}(\GrA)$. 
Thus, we have 
\begin{equation*}
\cD_{[\geq w]}(\Gr(\cA))|_U \, \subset \, \cD_{[\geq w]}(\GrA)
\end{equation*}

By Lemma \ref{Dlw_nonaffine}, we also have
\begin{equation*}
\cD_{< w]}(\Gr(\cA))|_U \, \subset \, \cD_{< w}(\GrA)
\end{equation*}
Therefore the restriction of \eqref{weight_decomp_seq_nonaffine} to $U$ becomes precisely the first row of \eqref{weight_trunc_exact_tri}.
This allows us to verify properties of weight truncations locally.

Combining the weight truncation sequence \eqref{weight_decomp_seq_nonaffine} with the local cohomology sequence \eqref{RGam_Ce_seq_sheaf}, this allows us to extend \eqref{triple_SOD_terms} to the case of non-affine base:
\begin{equation}  \label{triple_SOD_terms_nonaffine}
\begin{tikzcd} [row sep = 12, column sep = 15]
\cL_{[\geq w]} \RGam_{\scI^+}(\cM)  \ar[rr, "\cL_{[\geq w]} (\epsilon_{\cM})"] 
& &  \cL_{[\geq w]} \cM   \ar[rr, "\text{counit}"] 
\ar[ld, " \cL_{[\geq w]}(\eta_\cM)"]
& & \cM \ar[ld, "\text{unit}"] \\
& \cL_{[\geq w]} \Ce_{\scI^+}(\cM)  \ar[ul, "\text{[1]}" description]
& & \cL_{<w}(\cM) \ar[ul, "\text{[1]}" description]
\end{tikzcd}
\end{equation}
which gives the following generalization of \eqref{SOD_no_F} to the non-affine case:
\begin{equation*}
\cD(\Gr(\cA)) \, = \, \langle \, \cD_{<w}(\Gr(\cA)) \, , \, \cL_{[\geq w]} (\cD_{\cIpTR}(\Gr(\cA))) \, , \,  \cD_{\Torp,[\geq w]}(\Gr(\cA))  \, \rangle
\end{equation*}
where the component in the middle is the essential image of the functor $\cL_{[\geq w]} : \cD_{\cIpTR}(\Gr(\cA)) \ra \cD(\Gr(\cA))$, which is fully faithful with left quasi-inverse $\Ce_{\scI^+}$.
Alternatively, it may be described as
\begin{equation*}
\cL_{[\geq w]} (\cD_{\cIpTR}(\Gr(\cA))) \, = \, 
\{ \, \cM \in \cD_{[\geq w]}(\Gr(\cA)) \, | \, \RGam_{\scI^+}(\cM) \in \cD_{<w}(\Gr(\cA)) \, \}
\end{equation*}

Again, once we have formally defined the relevant functors, their properties can be checked affine-locally. For example, the analogues of Theorem \ref{three_term_SOD_Dmcoh}, \ref{I_triv_weight_equiv_Dsuitcoh}, \ref{cond_then_reg}, and \ref{smooth_then_reg_and_perf}  hold without change in the setting \eqref{sheaf_A_setting}.

\brm
As we explained in the introduction (see the paragraph of \eqref{birat_cobord_log_flip_diag_intro}), 
every wall-crossing in birational cobordism can be reduced to the setting \eqref{sheaf_A_setting} of a sheaf of $\bZ$-graded ring. Moreover, the ``master space construction'' of Thaddeus \cite{Tha96} realizes every variation of GIT quotient as a birational cobordism.
\erm

\section{Derived categories under flips and flops}  \label{Dcat_flip_flop_sec}

Let $A$ be a Noetherian $\bZ$-graded ring.
Recall from Proposition \ref{coh_IpTR_QpGr} and Remark \ref{stacky_Proj} that the derived category $\Dbcoh(\PProj^+(A))$ of the stacky projective space is equivalent to $\cD^b_{\coh(\IpTR)}(\GrA)$ defined in Definition \ref{IpTR_coh_def}. 
Thus, Corollary \ref{weight_Cech_Dsuitcoh_SOD} and Theorem \ref{I_triv_weight_equiv_Dsuitcoh} combine to show that $\Dbcoh(\PProj^+(A))$ is a semi-orthogonal summand of $\Dbcoh(\Fgw)$.
By symmetry, all these results hold for the negative direction as well, where the weight truncation is controlled by the full subcategory $\Flmw \subset \cF$. Thus, we see likewise that $\Dbcoh(\PProj^-(A))$ is a semi-orthogonal summand of $\Dbcoh(\Flmw)$.

As we observed in \eqref{transpose_Fga}, the pre-additive category $\Flmw$ is isomorphic to the opposite of $\Fgw$. Thus, $\Dbcoh(\PProj^+(A))$ is a semi-orthogonal summand of the derived category $\Dbcoh(\Fgw)$ of \emph{right} modules; while $\Dbcoh(\PProj^-(A))$ is a semi-orthogonal summand of the derived category $\Dbcoh((\Fgw)^{\op})$ of \emph{left} modules.
This suggests one to relate the derived categories by taking a duality functor (see \eqref{derived_Hom_ten_1} and \eqref{transpose_Fga})
\begin{equation}  \label{D_Fgw_functor}
\bD_{\Fgw} \, : \, \cD(\Fgw)^{\op} \raq \cD(\Flmw)\, , \qquad \cM \, \mapsto \, \bD_{\Fgw}(\cM) := \RHom_{\Fgw}(\cM, \Fgw)^{\tau}
\end{equation}

We will see that this tends to works well when there is a certain duality between the local cohomology complexes $\RGam_{I^+}(A)$ and $\RGam_{I^-}(A)$. Let $\omega_Y^{\bullet} \in \Dbcoh(A_0)$ be a dualizing complex, and take the weight degreewise dualizing functor
\begin{equation*}
\bD_Y : \cD(\GrA)^{\op} \ra \cD(\GrA) \, , \qquad \quad \bD_Y(M)_i \simeq \RHom_{A_0}(M_{-i}, \omega_Y^{\bullet})
\end{equation*} 
Then we will consider the assumptions
\begin{equation}  \label{Gor_and_local_cohom_dual_2}
\parbox{40em}{(i) $A$ is Gorenstein. \\
	(ii) There is an isomorphism $\Psi : \RGam_{I^+}(A)(a)[1] \xra{\cong }\bD_Y(\RGam_{I^-}(A))$ in $\cD(\GrA)$.}
\end{equation}
In fact, one of the main results in \cite{Yeu20a} is that a large class of flips and flops are controlled by sheaves of rings that satisfy this assumption, where $a = 0$ for a flop, and $a = 1$ for a flip.
We will see below that $\cD_{\IpTR}(\GrA)$ tends to be smaller than $\cD_{\ImTR}(\GrA)$ if $a > 0$, and they tend to be equivalent if $a = 0$. 
For some other results along this line, see \cite[Section 6]{Yeu20a}.

First, we observe that the duality functor \eqref{D_Fgw_functor} can be expressed alternatively in terms of a duality on $\cD(\GrA)$. Indeed, consider the functor
\begin{equation*}
\DAO \, : \, \cD(\GrA)^{\op} \raq \cD(\GrA) \, , \qquad \DAO(M) := \RHomcom_A(M,A) 
\end{equation*}
Then the following diagram commutes up to isomorphism of functors
\begin{equation}  \label{DAO_DFgw_diag}
\begin{tikzcd}[column sep = 80]
\cD_{[\geq w]}(\GrA)^{
	\op} \ar[r, " \cL_{[\leq -w]} \circ \DAO "] \ar[d, "(-)^{\sharp}_{\geq w}"', "\simeq"] & \cD_{[\leq -w]}(\GrA) \ar[d, "(-)^{\sharp}_{\leq -w}", "\simeq"']\\
\cD(\Fgw)^{\op} \ar[r, "\bD_{\Fgw}"] &  \cD(\cF_{[\leq -w]})
\end{tikzcd}
\end{equation}
Indeed, it suffices to prove the commutativity after replacing the vertical arrow on the left by its inverse $\scL_{[\geq w]} : \cD(\Fgw) \xra{\simeq} \cD_{[\geq w]}(\GrA)$, so that the commutativity follows (see \eqref{cF_A_tensor_Hom})
from taking the isomorphism
\begin{equation*}
\RHom_{\cF}(\cM \otimes_{\Fgw}^{{\bm L}} \cF  ,\cF) \, \cong \, \RHom_{\Fgw}(\cM, \cF)
\end{equation*}
in $\cD((\cF)^{\op})$, restrict to $\cD((\Fgw)^{\op})$, and take the transpose \eqref{transpose_Fga}. 

\bpp  \label{D_Fgw_ITR}
If \eqref{Gor_and_local_cohom_dual_2}(ii) holds for $a \geq 0$, then $\bD_{\Fgw}$ sends $\cD^-_{\coh,\IpTR}(\Fgw)$ to $\cD^+_{\coh,\ImTR}(\Flmw)$. 

If \eqref{Gor_and_local_cohom_dual_2}(ii) holds for $a = 0$, then $\bD_{\Flmw}$ also sends $\cD^-_{\coh,\ImTR}(\Flmw)$ to $\cD^+_{\coh,\IpTR}(\Fgw)$.
\epp

\bpf
It suffices to prove the first statement as the second statement follows by symmetry. Recall that the equivalence $(-)^{\sharp}_{\geq w} : \cD_{[\geq w]}(\GrA) \xra{\simeq} \cD(\Fgw)$ restricts to functors $\Dsuit_{\coh,[\geq w]}(\GrA) \ra \Dsuit_{\coh}(\Fgw)$ for each $\spadesuit \in \{ \, \, , +,-,b\}$, which is moreover an equivalence if $\spadesuit = -$.
Recall also that $M^{\sharp}_{\geq w} \cong (\cL_{[\geq w]} M)^{\sharp}_{\geq w}$ for all $M \in \cD(\GrA)$. Moreover, we have $M^{\sharp}_{\geq w} \in \cD_{\IpTR}(\Fgw)$ if and only if $\RGam_{I^+}(M) \in \cD_{<w}(\GrA)$.
Thus, in view of \eqref{DAO_DFgw_diag}, the result follows from the following Lemma.
\epf

\blm  \label{RGam_w_DAO}
Suppose that \eqref{Gor_and_local_cohom_dual_2}(ii) holds for $a \geq 0$. If $M \in \Dmcoh(\GrA)$ satisfies $\RGam_{I^+}(M) \in \cD_{<w}(\GrA)$, then we have
$\RGam_{I^-}(\DAO(M)) \in \cD_{>-w+a}(\GrA)$.
\elm

\bpf
Since $\bD_Y$ is involutive on complexes with locally coherent cohomology (see Lemma \ref{RGam_Dbcoh}), condition \eqref{Gor_and_local_cohom_dual_2}(ii) can be rewritten as an isomorphism $\RGam_{I^-}(A) \cong \bD_Y( \RGam_{I^+}(A) )(-a)[-1]$.
As a result, we have the following isomorphism in $\cD(\GrA)$:
\begin{equation*}
\begin{split}
\RHomcom_A( M , A ) \otimes_A^{{\bm L}} \RGam_{I^-}(A) 
\, &\cong \, \RHomcom_A( M , \RGam_{I^-}(A) ) \\
&\cong \, \RHomcom_A( M , \bD_Y( \RGam_{I^+}(A) ))(-a)[-1] \\
& = \, \bD_Y( \, \RGam_{I^+}(M) \,)(-a)[-1]
\end{split}
\end{equation*}
where the first isomorphism uses Proposition \ref{tensor_in_Hom_target}.
\epf

In Proposition \ref{D_Fgw_ITR}, if the small pre-additive category $\Fgw$ is ``Gorenstein'' in a suitable sense, then the functor $\bD_{\Fgw}$ gives a contravariant equivalence between $\Dbcoh(\Fgw)$ and $\Dbcoh(\Flmw)$. Proposition \ref{D_Fgw_ITR} then shows that $\cD^b_{\coh(\IpTR)}(\GrA)$ and $\cD^b_{\coh(\ImTR)}(\GrA)$ are related in the expected way: \ie,  the former is smaller for flips, and the two are equivalent for flops.

While $\Fgw$ may not be ``Goresntein'' in general, the duality functor $\bD_{\Fgw}$ is always well-behaved on $\Dperf(\Fgw)$. As a result, we have the following

\bcor  \label{flip_flop_Dperf}
Denote by $\cD_c \subset \cD$ the subcategory of compact objects.

Suppose that
\eqref{Gor_and_local_cohom_dual_2}(ii) holds for $a \geq 0$.
If the semi-orthogonal decomposition \eqref{Dbcoh_Fgw_SOD} on $\Dbcoh(\Fgw)$ restricts to one on $\Dperf(\Fgw)$, 
then the functor \eqref{D_Fgw_functor} restricts to a fully faithful functor
$\bD_{\Fgw} : (\cD_{\IpTR}(\Fgw)_c)^{\op} \rinto \cD_{\ImTR}(\Flmw)_c$.

Suppose that \eqref{Gor_and_local_cohom_dual_2}(ii) holds for $a = 0$. If the semi-orthogonal decomposition \eqref{Dbcoh_Fgw_SOD}, and its negative version, restrict to $\Dperf(\Fgw)$ and $\Dperf(\Flmw)$ respectively, then the functor \eqref{D_Fgw_functor} restricts to an equivalence
$\bD_{\Fgw} : (\cD_{\IpTR}(\Fgw)_c)^{\op} \xra{\simeq} \cD_{\ImTR}(\Flmw)_c$.
\ecor

\brm  \label{pairing_remark}
In Corollary \ref{flip_flop_Dperf}, we required \eqref{Dbcoh_Fgw_SOD} to restrict to $\Dperf(\Fgw)$ because the duality functor $\bD_{\Fgw}$ is well-behaved there. There is an alternative formulation of duality in terms of pairings of DG categories, which is well-behaved on the entire category. Namely, let $\Omega : \cA \times \cB \ra \cD(k)$ be a bi-cocontinuous exact functors between compactly generated triangulated category, with given DG enhancements. Then $\Omega$ induces co-continuous functors $\Omega^L : \cA \ra \cB^{\vee}$ and $\Omega^R : \cB \ra \cA^{\vee}$, where $\cB^{\vee}$ is the Ind-completion of $(\cB_c)^{\op}$, carried out at the DG level, and similarly for $\cA$. One can show that, if the essential image of $\Omega^L$ (resp. $\Omega^R$) contains all the compact objects, then $\Omega^R$ (resp. $\Omega^L$) is fully faithful. 
To apply this general formalism to our situation, consider $C^{(i)}_{\pm} \in \Dmcoh(\GrA)$ defined by
\begin{equation*}
C^{(i)}_+ := \cL_{[\geq w]}(\Ce_{I^+}(A)(-i)) 
\qquad \text{and} \qquad 
C^{(i)}_- := \cL_{[\leq -w]}(\Ce_{I^-}(A)(-i)) 
\end{equation*} 
One can show that, if the canonical evaluation maps
\begin{equation}  \label{evaluation_map_weight_trunc}
\begin{split}
\RHomcom_A( C^{(i)}_+ , A) \otimes_A^{{\bm L}} C^{(j)}_+ &\raq \RHomcom_A( C^{(i)}_+ , C^{(j)}_+ ) \\
\RHomcom_A( C^{(i)}_- , A) \otimes_A^{{\bm L}} C^{(j)}_- &\raq \RHomcom_A( C^{(i)}_- , C^{(j)}_- )
\end{split}
\end{equation}
are quasi-isomorphism in weight $0$, then we have
\begin{equation*}
\parbox{42em}{(1) If \eqref{Gor_and_local_cohom_dual_2}(ii) holds for $a \geq 0$, there is a fully faithful $\cD_{\IpTR}(\GrA) \rinto \cD_{\ImTR}(\GrA)$. \\
(2) If \eqref{Gor_and_local_cohom_dual_2}(ii) holds for $a = 0$, there is an equivalence $\cD_{\IpTR}(\GrA) \simeq \cD_{\ImTR}(\GrA)$.}
\end{equation*}
If the semi-orthogonal decomposition \eqref{Dbcoh_Fgw_SOD}, and its negative version, restrict to $\Dperf(\Fgw)$ and $\Dperf(\Flmw)$ respectively, then $C^{(i)}_{\pm}$ are in $\Dperf(\GrA)$, so that \eqref{evaluation_map_weight_trunc} are quasi-isomorphisms. In general, the evaluation map \eqref{evaluation_map_weight_trunc} being a quasi-isomorphism is essentially an issue about convergence of spectral sequences, and seems to be similar to the issues one encounters in Koszul duality. It seems possible that a formal modification of our arguments might lead to a more satisfactory statement. This might open up a way to tackle a conjecture of Bondal and Orlov.
\erm

\appendix

\section{Modules over pre-additive categories}  \label{app_mod_preadd}

A \emph{pre-additive category} is a category $\cA$ enriched over the monoidal category $(\Ab, \otimes)$ of abelian groups. It is said to be \emph{small} if the objects of $\cA$ form a set $\Ob(\cA)$.
It is helpful to think of a small pre-additive category as an ``associative ring with many objects'', as in \cite{Mit72}. This allows us to define the notions of left/right modules, tensor products, Hom spaces, etc, which we recall now.

Given a small pre-additive category $\cA$, a \emph{left $\cA$-module} is an additive functor $\cA \ra \Ab$, while a \emph{right $\cA$-module} is an additive functor $\cA^{\op} \ra \Ab$. 
Maps between left or right modules are simply natural transformations.
We will mostly work with right modules, and we denote the category of right $\cA$-modules by $\Mod(\cA)$. In more concrete terms, a right $\cA$-module associates an abelian group $M_a$ to each $a \in \Ob(\cA)$, together with maps $M_a \otimes \cA(a',a) \ra M_{a'}$, satisfying the obvious associativity and unitality conditions.

Given small pre-additive categories $\cA$ and $\cB$, 
an \emph{$(\cA,\cB)$-bimodule} consists of a collection $M(b,a) = \!\,_a M_b$ of abelian groups, one for each pair $a \in \Ob(\cA)$ and $b \in \Ob(\cB)$, together with maps $\cA(a,a') \otimes \!\,_a M_b \otimes \cB(b',b) \ra \!\,_{a'} M_{b'}$, satisfying the obvious associativity and unitality conditions. For example, $\cA$ is canonically a bimodule over itself.
Denote by $\!\,_{\cA}\Mod_{\cB}$ the category of $(\cA,\cB)$-bimodules.

If $M \in \!\,_{\cA}\Mod_{\cB}$ and $N \in \!\,_{\cB}\Mod_{\cC}$, then define
$M \otimes_{\cB} N \in \!\,_{\cA}\Mod_{\cC}$ by
\begin{equation}  \label{bimod_tensor}
\! \,_a(M \otimes_{\cB} N)_{c} \, := \bigl( \, \bigoplus_{b \in \Ob(\cB)} \! \,_aM_b \otimes \! \,_b N_c \bigr) \big/ ( \, \xi f \otimes \eta - \xi \otimes f \eta \,)
\end{equation} 
where we mod out the abelian subgroup generated by the displayed relations, for $\xi \in \,_aM_{b'}$, $f \in \cB(b,b')$, and $\eta \in \,_bM_{c}$.
In particular, if $\cA = \cC = \ast$ is the pre-additive category with one object, with endomorphism algebra $\bZ$, then this gives the notion of a tensor product $M \otimes_{\cB} N \in \Ab$ between a right $\cB$-module $M$ and a left $\cB$-module $N$.

Similarly, if $M \in \!\,_{\cA}\Mod_{\cB}$  and $N \in \!\,_{\cC}\Mod_{\cB}$, then we define  $\Hom_{\cB}(M,N) \in \!\,_{\cC}\Mod_{\cA}$ by
\begin{equation}   \label{bimod_Hom}
\! \,_c \Hom_{\cB}(M,N)_{a} \, := \, \bigl \{ (\varphi_b) \in \prod_{b \in \Ob(\cB)} \Hom_{\Ab}( \! \,_a M_b,  \! \,_c N_b) \, \big| \, \varphi_{b}(\xi f) = \varphi_{b'} (\xi) f  \, \, \,  \forall \, 
\xi \in \! \,_a M_{b'} , \,  f \in \cB(b,b') \,\bigr \}
\end{equation} 
In other words, $\! \,_c \Hom_{\cB}(M,N)_{a}$ is the Hom-space $\Hom_{\cB}(\!\,_a M , \!\,_c N)$ in the (big) additive category $\Mod(\cB)$.

As for usual associative algebras, there are canonical isomorphisms
\begin{equation} \label{triv_Hom_ten}
M \otimes_{\cB} \cB  \cong   M   \cong  \cA \otimes_{\cA} M  
\qquad \text{ and } \qquad 
\Hom_{\cB}(\cB,M) \cong M
\end{equation}

For any $M \in \!\,_{\cA}\Mod_{\cB}$, $N \in \!\,_{\cB}\Mod_{\cC}$, and $L \in \!\,_{\cE}\Mod_{\cC}$, there is a usual Hom-tensor adjunction, given by the canonical isomorphism of $(\cE,\cA)$-bimodules
\begin{equation}  \label{bimod_Hom_ten_adj}
\Hom_{\cC}(M \otimes_{\cB} N , L) \, \cong \, \Hom_{\cB}(M , \Hom_{\cC}(N,L))
\end{equation}

For each $a \in \Ob(\cA)$, denote by $\!\,_a\cA$ the right $\cA$-module represented by $a$. In other words, $\!\,_a\cA_{a'} := \cA(a',a)$. A right module $M$ is said to be \emph{free} if there is an indexed set of objects $\varphi : S \ra \Ob(\cA)$, together with an isomorphism $M \cong \oplus_{s \in S}  \, (\!\,_{\varphi(s)}\cA)$.
In more concrete terms, this means that there is a set $S$ of elements $\xi_s \in M_{\varphi(s)}$ such that, for any $a \in \Ob(\cA)$, any element $\xi \in M_a$ can be written uniquely as a finite sum $\xi = \sum \xi_s f_s$, for $f_s \in \cA(a,\varphi(s))$.
The cardinality of $S$ is said to be the \emph{rank} of the free module $M$.

Clearly the category $\Mod(\cA)$ of right modules is an abelian category, where limits and colimits are determined objectwise. Thus, it also satisfies the usual Ab5 and Ab3* axioms of an abelian category. Moreover, the set $\{ \,_a\cA \, \}_{a \in \Ob(\cA)}$ of right modules forms a set of generators for $\Mod(\cA)$, so that $\Mod(\cA)$ is a Grothendieck category (see, e.g., \cite[Tag 079B]{Sta}). 
Projective objects in $\Mod(\cA)$ are precisely retracts of free modules. A projective right module is said to be of \emph{finite rank} if it is a retract of a free module of finite rank.

\bdf \label{fin_gen_mod}
We say that a right module $M \in \Mod(\cA)$ is \emph{finitely generated} if there is an epimorphism 
$\oplus_{s \in S}  \, (\!\,_{\varphi(s)}\cA) \ronto M$ for a finite indexed set of objects $\varphi : S \ra \Ob(\cA)$.

In more concrete terms, this means that there is a finite set $S$ of elements $\xi_s \in M_{\varphi(s)}$ such that, for any $a \in \Ob(\cA)$, any element $\xi \in M_a$ can be written as a finite sum $\xi = \sum \xi_s f_s$, for $f_s \in \cA(a,\varphi(s))$.
\edf

\bdf
A small pre-additive category $\cA$ is said to be \emph{right Noetherian} (resp. \emph{left Noetherian}) if every submodule of a finitely generated right (resp. left) $\cA$-module is finitely generated.
It is said to be \emph{Noetherian} if it is both left and right Noetherian.
\edf

Since $\Mod(\cA)$ is a Grothendieck category, it has enough injectives (see, e.g., \cite[Tag 079H]{Sta}). Moreover, complexes in $\Mod(\cA)$ admit K-injective resolutions (see, e.g., \cite[Tag 079P]{Sta}). The category $\Mod(\cA)$ clearly has enough projectives. Thus, by 
\cite[Theorem 3.4]{Spa88} (see also \cite[Tag 06XX]{Sta}), complexes in $\Mod(\cA)$ admit K-projective resolutions. This allows us to take derived functors of the above Hom functors and tensor functors.
However, a subtlety arises when one wants to take the derived tensor product or the derived Hom bimodule between bimodules. For example, even if $M \in \!\,_{\cA}\Mod_{\cB}$ is projective in the category $\!\,_{\cA}\Mod_{\cB}$, it may not be true that each $\!\,_a M \in \Mod(\cB)$ is projective, or even flat, so that it might be problematic if one wants to derive \eqref{bimod_tensor} and \eqref{bimod_Hom} naively. However, for our purposes, we will only need to consider the derived tensor product (or derived Hom) between a module and a bimodule. For these, there are no problems, and we may define
\begin{equation}  \label{derived_Hom_ten_1}
\begin{split}
- \otimes_{\cA}^{{\bm L}} - \, &: \, \cD(\cA) \, \times \, \cD( \!\,_{\cA}\Mod_{\cB} ) \raq \cD(\cB) \\
\RHom_{\cA}(-,-) \, &: \, \cD(\cA)^{\op} \, \times \, \cD( \!\,_{\cB}\Mod_{\cA} ) \raq \cD(\cB^{\op})
\end{split}
\end{equation}

In particular, given an additive functor $F : \cA \ra \cB$, then $\cB$ may be viewed as an $(\cA,\cB)$-bimodule in the obvious way, so that the extension functor
\begin{equation*}
- \otimes_{\cA}^{{\bm L}} \cB \, : \, \cD(\cA) \raq \cD(\cB) 
\end{equation*}
is well-defined. Moreover, it satisfies the usual ``cancellation rule'' for tensor products
\begin{equation}  \label{tensor_two_step}
(M \otimes^{{\bm L}}_{\cA} \cB) \otimes^{{\bm L}}_{\cB} N \, \cong \, M \otimes^{{\bm L}}_{\cA} N
\end{equation}
for any $M \in \cD(\cA)$ and $N \in \cD(\!\,_{\cB}\Mod_{\cC})$.

\bdf
An object $M \in \cD(\cA)$ is said to be \emph{pseudo-coherent} if it is quasi-isomorphic to a bounded above complex $P^{\bullet}$ of projective modules of finite rank. 
Denote by $\Dpc(\cA) \subset \cD(\cA)$ the full subcategory consisting of pseudo-coherent objects.
\edf

\bdf
Suppose $\cA$ is right Noetherian, then denote by $\Dmcoh(\cA) \subset \cD(\cA)$ the full subcategory consisting of objects $M \in \cD(\cA)$ such that each $H^p(M)$ is finitely generated, and $H^p(M) = 0$ for $p \gg 0$.
\edf

\bpp   \label{Dmcoh_Dpc_right_Noeth}
Suppose $\cA$ is right Noetherian, then for any $M \in \cD(\cA)$, the followings are equivalent:
\begin{enumerate}
	\item $M \in \Dpc(\cA)$;
	\item $M \in \Dmcoh(\cA)$;
	\item $M$ is quasi-isomorphic to a bounded above complex of free modules of finite rank.
\end{enumerate}
\epp

\bpf
Clearly only the implication $(2) \Rightarrow (3)$ needs proof. It follows from the well-known Lemma \ref{Dbcoh_pc_lem} below (for a proof, see, \eg, \cite[Lemma A.40]{Yeu20c}).
\epf

\blm  \label{Dbcoh_pc_lem}
Let $\cC$ be an abelian category, and let $\cP \subset \Ob(\cC)$ be a collection of projective objects closed under finite direct sum. Denote by $Q(\cP) \subset \Ob(\cC)$ the collection of objects $M$ such that there exists an epimorphism $P \ronto M$ from some $P \in \cP$. 
Suppose $Q(\cP)$ is closed under taking subobjects, then for any bounded above complex $M^{\bullet}$ in $\cC$ whose cohomology objects lie in $Q(\cP)$, there exists a bounded above complex $P^{\bullet}$ of objects in $\cP$, together with a quasi-isomorphism $\varphi: P^{\bullet} \rsa M^{\bullet}$.
\elm

Clearly, the set $\{ \!\,_a \cA \}_{a \in \Ob(\cA)}$ forms a set of compact generators of $\cD(\cA)$. Denote by $\Dperf(\cA)$ the smallest split-closed triangulated subcategory of $\cD(\cA)$ containing the set $\{ \!\,_a \cA \}_{a \in \Ob(\cA)}$ of objects, then it is a standard fact (see, e.g., \cite[Theorem 4.22]{Rou08} or \cite[Lemma 2.2]{Nee92}) that 
\begin{equation}  \label{cpt_equal_perf}
\cD(\cA)_c \, = \, \Dperf(\cA)
\end{equation}
where the subscript $(-)_c$ denotes the subcategory of compact objects.

\end{document}